\documentclass[a4paper,leqno]{amsart}

\usepackage{amsmath}
\usepackage{amsthm}
\usepackage{amssymb}
\usepackage{amsfonts}
\usepackage[applemac]{inputenc}
\usepackage{mathrsfs}
 \usepackage{pdfsync}
\usepackage{xcolor}
\usepackage{graphicx,epstopdf,color}

\def\R{{\mathbb R}}
\def\N{{\mathbb N}}
\def\({\left(}
\def\){\right)}
\def\<{\left\langle}
\def\>{\right\rangle}
\def\eps{\varepsilon}

\def\le{\leqslant}
\def\ge{\geqslant}

\DeclareMathOperator{\IM}{Im}

\def\Tend#1#2{\mathop{\longrightarrow}\limits_{#1\rightarrow#2}}

\theoremstyle{plain}
\newtheorem{theorem}{Theorem}[section]

\newtheorem{proposition}[theorem]{Proposition}

\newtheorem{conjecture}[theorem]{Conjecture}

\theoremstyle{definition}
\newtheorem{definition}[theorem]{Definition}

\newtheorem{remark}[theorem]{Remark}
\newtheorem*{remark*}{Remark}

\numberwithin{equation}{section}

\begin{document}

\title[Cubic-quintic ground states]
{On ground state (in-)stability in multi-dimensional
cubic-quintic Schr\"odinger equations}

\author[R. Carles]{R\'emi Carles}
\address[R. Carles]
{Univ Rennes, CNRS\\ IRMAR - UMR 6625\\ F-35000
  Rennes, France}
\email{Remi.Carles@math.cnrs.fr}

\author[C. Klein]{Christian Klein}
\address[C.~Klein]
{Institut de Math\'ematiques de Bourgogne, Universit\'e de 
Bourgogne-Franche-Comt\'e, 
9 avenue Alain Savary, BP 47870, 21078 Dijon Cedex}
\email{christian.klein@u-bourgogne.fr}

\author[C. Sparber]{Christof Sparber}
\address[C.~Sparber]
{Department of Mathematics, Statistics, and Computer Science, M/C 249, University of Illinois at Chicago, 851 S. Morgan Street, Chicago, IL 60607, USA}
\email{sparber@math.uic.edu}

\begin{abstract}
We consider the nonlinear Schr\"odinger equation with a focusing cubic
term and a defocusing quintic nonlinearity in dimensions two and
three. The main interest of this article is the problem of orbital (in-)stability of ground state solitary
waves. We recall the notions of energy minimizing versus action minimizing ground states 
and prove that, in general, the two must be considered as nonequivalent. We numerically investigate the 
orbital stability of least action ground states in the radially 
symmetric case, confirming existing conjectures or
leading to new ones. 
\end{abstract}

\date{\today}


\thanks{This publication is based on work supported by the NSF 
  through grant no. DMS-1348092.
RC acknowledges support from Rennes M\'etropole, through its
  AIS program.
CK is partially supported by 
the ANR-FWF project ANuI - ANR-17-CE40-0035, the isite BFC project 
NAANoD, the EIPHI Graduate School (contract ANR-17-EURE-0002), by the 
European Union Horizon 2020 research and innovation program under the 
Marie Sklodowska-Curie RISE 2017 grant agreement no. 778010 IPaDEGAN 
and the EITAG project funded by the FEDER de Bourgogne, the region 
Bourgogne-Franche-Comt\'e and the EUR EIPHI}

\maketitle


\section{Introduction}
\label{sec:intro}

\subsection{Basic setting} This work is concerned with the time-evolution corresponding to the {\it cubic-quintic nonlinear Schr\"odinger equation} (NLS)
\begin{equation}
  \label{eq:nls}
  i\partial_t u +\frac{1}{2}\Delta u =-|u|^2u+|u|^4u,\quad (t, x)\in \R\times \R^{d},
\end{equation}
in dimensions $d=2$, or $d=3$, and subject to initial data 
\[u_{\mid
  t=0}=u_0\in H^1(\R^d).
  \]  
  
The quintic modification of the cubic
Schr\"odinger equation is a model which was introduced in the
one-dimensional case in \cite{Pushkarov}, as an approximate model in
the framework of nonlinear optics. Equation~\eqref{eq:nls}  appeared
more recently in the context of Bose--Einstein condensation, with
$d=2$ or $3$: see
e.g. \cite{PhysRevA63,JPhysB,PhysRevE}, and \cite{Malomed19} for a
review. 
In space dimensions $d=2$ or
$3$, the impact of the quintic term on the dynamical properties of the solution $u$ 
is stronger than in $d=1$, as we shall discuss below. 

Depending on the space dimension, which we always assume at most three
to simplify the discussion, the nonlinearity in this model is
seen to be: focusing $L^2$-subcritical plus defocusing
$L^2$-critical ($d=1$), 
focusing $L^2$-critical plus defocusing
$H^1$-subcritical ($d=2$), or focusing $L^2$-supercritical plus
defocusing $H^1$-critical ($d=3$). Recall that for the purely focusing cubic
NLS, solitons exist in every dimension and {\it finite time blow-up} is possible provided $d\ge 2$
(see e.g. \cite{CazCourant}). The presence of the quintic nonlinearity
prevents finite time blow-up in $d=2$ or $3$ (see
Proposition~\ref{prop:GWP} below), and also affects the stability of
solitary waves: understanding this latter aspect more precisely is the
main motivation for this paper. For a more precise discussion on the role of criticality in combined power 
nonlinearities see \cite{TaoVisanZhang}.

\smallbreak

The NLS \eqref{eq:nls} formally enjoys the following basic conservation laws:

\begin{enumerate}
\item Mass: \[
  M(u)=\| u(t, \cdot) \|_{L^2(\R^d)}^2,\]

\item  Momentum: \[
 P(u)=\IM\int_{\R^d}\bar u(t,x)\nabla u(t,x)dx,\]

\item Energy:
\[
 E(u) = \frac{1}{2}\|\nabla u(t, \cdot)\|_{L^2(\R^d)}^2 -\frac{1}{2}\|
u(t, \cdot)\|_{L^4(\R^d)}^4+ \frac{1}{3}\|  u(t, \cdot)\|_{L^6(\R^d)}^6.\]
\end{enumerate}
As evoked above, one important effect of the defocusing, quintic term
is to {\it prevent finite time blow-up} which may occur in the purely cubic
case. Indeed, the conservation of the energy, combined with H\"older's inequality,
\begin{equation}\label{eq:holder}
  \|u \|_{L^4(\R^d)}^4 \le \|u \|_{L^2(\R^d)}\|u \|_{L^6(\R^d)}^3, 
\end{equation}
shows that the focusing, cubic part cannot be an obstruction to the
existence of a global in-time solution. More precisely, we have, in
view of \cite{CazCourant} for $d=2$ and \cite{Zhang06} for $d=3$:
\begin{proposition}[Global well-posedness]\label{prop:GWP}
  Let $d= 2, 3$. For any initial data $u_0\in H^1(\R^d)$, the equation \eqref{eq:nls} has a
  unique solution $u\in C(\R;H^1(\R^d))$, such that $u_{\mid
    t=0}=u_0$. This solution obeys the conservation of mass, energy, and
  momentum. 
\end{proposition}

We note that in
  \cite{LeMesurier1988}, numerical simulations are presented, in which the influence of a small defocusing quintic term on the 
  time-evolution of a focusing cubic NLS is studied. In $d=2$ and $3$, and for initial data consisting of Gaussians, one obtains 
  a time-periodic (multi-focusing) solution, similar to the one depicted in Fig. \ref{NL35_d2solom005}.  

\subsection{Orbital stability of action minimizing ground states}  A particular class of global solutions are time-periodic {\it solitary
waves} of the form $u(t,x)=e^{i\omega t}\phi(x)$, with $\omega\in \R$ and
  $\phi$ satisfying
  \begin{equation}\label{eq:soliton}
    -\frac{1}{2}\Delta \phi + \omega \phi
    -|\phi|^2\phi+|\phi|^4\phi=0,\quad \phi\in
    H^1(\R^d)\setminus\{0\}. 
  \end{equation}
For $d\le 3$, solitary waves $\phi$ exist provided that the frequency $\omega$
satisfies the (necessary and sufficient) condition:
\begin{equation*}
  0<\omega<\tfrac{3}{16},
\end{equation*}
see \cite{CaSp}. Given an admissible $\omega \in (0, \tfrac{3}{16})$, we may then look for {\it least action
ground states}, i.e. solutions $\phi_\omega(x)$ which minimize the {\it action} 
\[S_\omega(\phi) = E(\phi) + \omega M(\phi)\]
  among all nontrivial stationary solutions $\phi\in
  H^1(\R^d)$. Indeed, it is known from \cite{ByJeMa09,CiJeSe09} that
  every minimizer of the action  $S_\omega$ is of the  
form
\begin{equation}\label{eq:phiomega}
\phi_\omega(x) = e^{i\theta}Q_\omega(x-x_0),
\end{equation}
for some constant $\theta\in \R$, $x_0\in \R^d$, and with $Q_\omega$
the unique {\it positive, radial solution} to \eqref{eq:soliton}. In the following, we 
are mainly interested in the 
orbital stability of these specific solutions. To this end, we note that, as in the case of more standard, homogeneous nonlinearities
(e.g. the cubic case), the NLS \eqref{eq:nls} enjoys three important invariances:
\begin{itemize}
\item[(i)] Spatial translation: if $u(t,x)$ solves \eqref{eq:nls}, then so
  does $u(t,x-x_0)$, for any given $x_0\in \R^d$.
\item[(ii)] Gauge: if $u(t,x)$ solves \eqref{eq:nls}, then so
  does $e^{i\theta}u(t,x)$, for any given constant $\theta\in \R$.
\item[(iii)] Galilean: if $u(t,x)$ solves \eqref{eq:nls}, then so
  does $ u(t,x-vt)e^{iv\cdot x-i|v|^2 t/2}$, for any given $v\in
  \R^d$. 
\end{itemize}
The first two invariants are seen to be present in formula
\eqref{eq:phiomega}. In combination with the third one, these invariants motivate the following
standard notion of stability (see e.g. \cite{CazCourant}):
\begin{definition}\label{def:stability}
  Let $\phi$ be a solution of \eqref{eq:soliton}. The standing wave
  $e^{i\omega t}\phi(x)$ is {\it orbitally stable} in $H^1(\R^d)$, if
  for all $\eps>0$, there exists $\delta>0$ such that if $u_0\in
  H^1(\R^d)$ satisfies
  \[\|u_0-\phi\|_{H^1(\R^d)}\le \delta,\]
  then the
  solution to \eqref{eq:nls} with $u_{\mid t=0}=u_0$ satisfies
  \begin{equation*}
    \sup_{t\in \R}\inf_{{\theta\in \R}\atop{y\in
      \R^d}}\left\|u(t, \cdot)-e^{i\theta}\phi(\cdot
      -y)\right\|_{H^1(\R^d)}\le \eps.
  \end{equation*}
  Otherwise, the standing wave is said to be unstable. 
\end{definition}
Following the breakthrough due to
M.~Weinstein, Grillakis, Shatah and Strauss
introduced a general stability/instability criterion in \cite{GSS87} (see
also \cite{BGR15}). Assuming certain spectral properties of the linearization of
\eqref{eq:soliton} about $Q_\omega$ (which are satisfied in the
present cubic-quintic case, see
e.g. \cite{CaSp,KOPV17,LewinRotaNodari2015}), one has the following dichotomy:
\begin{itemize}
\item[(i)] If $\frac{\partial }{\partial \omega}M(Q_\omega)>0$, then
    $e^{i\omega t}Q_\omega(x)$ is orbitally stable,
    
 \item[(ii)] If $\frac{\partial}{\partial \omega}M(Q_\omega)<0$, then
    $e^{i\omega t}Q_\omega(x)$ is unstable.
\end{itemize}
This criterion has proven extremely useful in  the case of homogeneous
nonlinearities, as well as in the case of mixed nonlinearities in $d=1$ thanks to an explicit formula, cf. \cite{IlievKirchev93,Ohta95}. 
In particular, when $d=1$, all ground
states solitary waves for \eqref{eq:nls} are orbitally stable. However,
in the case $d=2$ or $3$, only partial results are currently available by using the criterion above, see Sections~\ref{sec:rappel2D}
and \ref{sec:rappel3D}, respectively. In our numerical simulations, we
will only consider radial 
  perturbation of ground states, and thus
  remain in the radial framework. The notion of orbital stability then
  coincides with asymptotic stability up to a phase.

\begin{remark}
  The dependence of $\delta$ upon $\eps$ in
  Definition~\ref{def:stability} is unknown, in general. The proof
  of stability via the above criterion provides a rather explicit
  dependence of (a possible) $\delta$ as a function of $\eps$, when
  $\tfrac{\partial }{\partial \omega}M(Q_\omega)>0$ is known. On the
  other hand, the stability of the set of constrained energy minimizers
  (cf. Definition~\ref{def:set-stability} below) is obtained by a
  non-constructive argument. In numerical
  simulations, tuning
  initial perturbations of a solitary wave which are sufficiently
  large to be visible, but not too large (to still adhere to the notion of
  stability), requires a subtle balance. 
\end{remark}


\subsection{Constrained energy minimizers} Since solitary waves may be obtained by other means than minimizing the action $S_\omega$, 
one may want to look for alternative approaches to 
orbital stability. An important such alternative is obtained if for $\rho>0$, we denote
  \begin{equation*}
    \Gamma(\rho) = \left\{ u\in H^1(\R^d),\  M(u)=\rho\right\},
  \end{equation*}
and assume that the constrained minimization problem
  \begin{equation}\label{eq:8.3.5}
u\in \Gamma(\rho),\quad E(u)=    \inf \{ E(v)\ ;\  v\in \Gamma(\rho)\}
=:E_{\rm min}(\rho) 
  \end{equation}
  has a solution. We call such minimizers
    \emph{energy ground states}, in order to make the distinction with
    least action ground states as clear as possible. Denote by $\mathcal
  E(\rho)$ the set of such 
  solutions, i.e.
  \begin{equation*}
    \mathcal  E(\rho):=\{u \in H^1(\R^d),\  M(u)=\rho,\ E(u) = E_{\rm
      min}(\rho)\}.
  \end{equation*}
  Now, let $\phi\in \mathcal E(\rho)$: Then there exists a Lagrange multiplier $\Lambda$ such
  that 
  \begin{equation*}
    dE(\phi) =\Lambda dM(\phi),
  \end{equation*}
  and thus, $\phi$ solves the stationary Schr\"odinger equation \eqref{eq:soliton} for some (unknown)
  $\omega\in (0, \tfrac{3}{16})$. Observe that if $\phi\in \mathcal E(\rho)$, then 
\[
\{ e^{i\theta}
\phi(\cdot -y);\ \theta\in\R,\ y\in \R^d\}\subset \mathcal
E(\rho).
\] 
When the nonlinearity is homogeneous (and
$L^2$-subcritical), this inclusion becomes an equality, see
\cite{CaLi82,CazCourant}. 
However, for non-homogeneous nonlinearities like in \eqref{eq:nls}, 
relating these two constructions of solitary waves 
  (i.e., action minimizing ground states versus constrained energy minimizers) is not
  obvious at all, and the issue is possibly more complex than it may
  appear at a glance. First, a-priori
  nothing guarantees that an element of $\mathcal 
E(\rho)$ minimizes the action. Second, and this is more subtle:  consider a least action ground
  state $Q_\omega$, and let 
  $\rho=M(Q_\omega)$. It is not obvious, and not necessarily true, that
  $Q_\omega\in \mathcal E(\rho)$. In particular, the map $\rho \mapsto \omega$ may not be 
  one-to-one. We will prove that, unlike in the case
  of homogeneous nonlinearities, we may indeed have
  $Q_\omega\not\in \mathcal E(\rho)$, cf. Theorem~\ref{theo:different} below. 
 This fact should be compared to the recent results of \cite{JeanjeanLu-p}: In there, the authors establish
  for a large class of nonlinearities (including the cubic-quintic
  one in space dimension $d\le 3$), that all \emph{all energy minimizing 
  ground states are least action ground states. In addition, they show that if $\omega$ is obtained 
  as the Lagrange multiplier associated to the mass constrained $M(u)=\rho$, then any least action solution of \eqref{eq:soliton}
  at this value of $\omega$ is a constrained energy minimizer with the same mass $\rho$.}

\smallbreak

This paper is now organized as follows: In Section \ref{sec:rev}, we shall review several known mathematical results on the 
(in-)stability of ground state solitary waves in $d=1,2,3$. In particular, we recall the fact that the set of 
energy minimizers is orbitally stable. We then prove that the dynamics of this set can be distinguished from 
dispersive behavior and that in $d=3$ the sets of action and energy minimizers are not equivalent.
In Section \ref{sec:numground} we numerically construct 
action  ground states and also collect several of their qualitative properties. Numerical evidence for the orbital stability of these action ground states in $d=2$
is then given in Section \ref{sec:stab2D}, where we will also describe the numerical algorithm used to simulate the time evolution 
of \eqref{eq:nls}. Finally, we shall turn to the question of orbital stability and instability of 3D action  ground states in Section \ref{sec:stab3D}, where 
we will provide numerical evidence for several conjectures on the particular nature of the instability.


\section{Mathematical results on orbital (in-)stability}\label{sec:rev}

\subsection{Stability for energy ground states} As a preliminary step, we shall recall the {\it Pohozaev identities} for the cubic-quintic
case (for a derivation, see e.g. \cite{CaSp}): if $\phi\in H^1(\R^d)$ solves the stationary Schrödinger equation \eqref{eq:soliton}, then
 \begin{equation}
  \label{eq:phi1}
  \frac{1}{2}\int_{\R^d}|\nabla \phi|^2 \, dx
  -\int_{\R^d}|\phi|^4 \, dx+\int_{\R^d}|\phi|^6 \, dx+\omega \int_{\R^d}|\phi|^2  \, dx= 0,
\end{equation}
as well as
\begin{equation}
  \label{eq:phi2}
  \frac{d-2}{2} \int_{\R^d}|\nabla \phi|^2 \, dx
 -\frac{d}{2}\int_{\R^d}|\phi|^4 \, dx+\frac{d}{3}\int_{\R^d}|\phi|^6 \, dx+\omega
 d \int_{\R^d}|\phi|^2  \, dx= 0.
\end{equation}
The aforementioned admissible range for $\omega \in (0, \tfrac{3}{16})$ is one of the consequences of these identities.
In addition, if $d=2$, and after multiplying \eqref{eq:phi1} by $2$
and subtracting \eqref{eq:phi2}, we find
\begin{equation*}
  0=\|\nabla\phi\|_{L^2(\R^2)}^2
  -\|\phi\|_{L^4(\R^2)}^4+\frac{4}{3}\|\phi\|_{L^(\R^2)}^6 = 2E(\phi)
  +\frac{5}{6}\|\phi\|_{L^(\R^2)}^6 . 
\end{equation*}
Therefore, any solitary wave in 2D has {\it negative energy}. In the 3D case,
this is not necessarily so, as we will see below.

Next, we recall the notion of orbital stability for the set of energy minimizers, as introduced in \cite{CaLi82}:

\begin{definition}\label{def:set-stability}
 We say that solitary waves are $\mathcal E(\rho)$-{\it orbitally
 stable}, if for all $\eps>0$, there exists $\delta>0$ such that if
 $u_0\in H^1(\R^d)$ satisfies
 \[\inf_{\phi\in \mathcal E(\rho)}\|u_0-\phi\|_{H^1(\R^d)}\le \delta,\]
  then the
  solution to \eqref{eq:nls} with $u_{\mid t=0}=u_0$ satisfies
  \begin{equation*}
    \sup_{t\in \R}\inf_{\phi\in \mathcal E(\rho)}\left\|u(t,
      \cdot)-\phi\right\|_{H^1(\R^d)}\le \eps. 
  \end{equation*}
\end{definition}
This notion is weaker than the one given in
Definition~\ref{def:stability}, in the sense that $\mathcal E(\rho)$
may be a large set. In particular, we cannot infer orbital stability of 
individual members $\phi\in \mathcal E(\rho)$, as required in Definition \ref{def:stability}.
In the case of homogeneous nonlinearities, however, 
one can prove that the set $\mathcal E(\rho)$ consists of only a
single element $\phi$ 
(up to translation and phase conjugation), and thus one recovers Definition \ref{def:stability} from the one above.
To prove $\mathcal E(\rho)$-orbital stability via the concentration-compactness method, the main step consists in showing that the minimal energy $E_{\rm min}(\rho)<0$. In
our case this yields:
\begin{theorem}[From \cite{CaSp}]\label{theo:stability-constrained}
  Let $d=2$ or $3$.
  \begin{enumerate}
    \item If $E_{\rm min}(\rho)<0$, then  $\mathcal
  E(\rho)$ is not empty, and
  the set of energy ground states is $\mathcal
  E(\rho)$-orbitally stable. 
  \item  There exists $\rho_0(d)>0$ such that for
    $\rho>\rho_0(d)$, $E_{\rm min}(\rho)<0$.
  \end{enumerate}
\end{theorem}
In particular, it seems reasonable to expect that no dispersion is possible near
elements of $\mathcal E(\rho)$, in the sense that a solution $u(t, \cdot)$ within 
Definition~\ref{def:set-stability} cannot satisfy
\begin{equation}\label{eq:disp-Linfty}
  \|u(t)\|_{L^\infty(\R^d)}\Tend t \infty 0. 
\end{equation}
We will use this criterion  as a guiding principle for interpreting several of our numerical findings  below. 
However, it is not clear \emph{a priori} that for $\rho>0$, such that
$\mathcal E(\rho)\not =\emptyset$, we have
\begin{equation*}
  \inf\{ \|\phi\|_{L^\infty(\R^d)},\ \phi\in \mathcal E(\rho)\}>0.
\end{equation*}
The property $E_{\rm min}(\rho)<0$ 
makes it possible to rule out this scenario. 
 
\begin{proposition}[Non-dispersion of energy ground
    states]\label{prop:ndisp}  
Let $d=2$ or $3$, and $\rho>0$. If $E_{\rm
    min}(\rho)<0$, then
\begin{equation*}
m_p:=  \inf\{ \|\phi\|_{L^p(\R^d)},\ \phi\in \mathcal E(\rho)\}>0,
\end{equation*}
for any $p\in (4, \infty]$. In particular,
there exists $\eps_0(d)$  such that for $0<\eps\le \eps_0(d)$, any
solution $u$ provided by  Definition~\ref{def:set-stability} satisfies
  \[
\inf_{t \in \R} \| u(t, \cdot) \|_{L^\infty(\R^d)} >0.
\]
\end{proposition}
The statement of this proposition does not involve the above parameter
$\rho_0(d)$. For practical application, we want to emphasize that if
for a given mass $\rho$, a stationary solution has negative energy,
then solutions around $\mathcal E(\rho)$ cannot disperse.
\begin{proof} 
Assume, by contradiction, that there exists a sequence $(\phi_n)_{n\in \N} \subset \mathcal E(\rho)$ such that 
\[
\lim_{n \to \infty} \| \phi_n\|_{L^p(\R^d)}=0, \quad \text{for some $p>4$.}
\]
Since $\| \phi_n\|_{L^2(\R^d)}=\rho$, by interpolation, this implies that 
\[
\lim_{n \to \infty} \| \phi_n\|_{L^4(\R^d)}=0.
\]
In turn, this yields that 
\[
E_{\rm min} (\rho)=\lim_{n \to\infty} E(\phi_n)= \lim_{n \to \infty} \(\frac{1}{2}\|\nabla \phi_n\|_{L^2(\R^d)}^2 + \frac{1}{3}\|  \phi_n\|_{L^6(\R^d)}^6\)\ge 0,
\]
a contradiction.

Now,  choose $0<\eps<m_5$. Consider initial data $u_0\in H^1(\R^d)$ such that 
 \[\inf_{\phi\in \mathcal E(\rho)}\|u_0-\phi\|_{H^1(\R^d)}\le \delta,\]
where $\delta$ stems from Definition~\ref{def:set-stability} and
orbital stability (which is ensured since $E_{\rm min}(\rho)<0$). 
 Then, by the $\mathcal E(\rho)$-orbital stability and Sobolev
 imbedding, we have  
\begin{equation*}
    \sup_{t\in \R}\inf_{\phi\in \mathcal E(\rho)}\left\|u(t,
      \cdot)-\phi\right\|_{L^5(\R^d)}\le \eps. 
  \end{equation*}
  Since for all $t$ and all $\phi\in \mathcal E(\rho)$,
  \begin{equation*}
    \| u(t, \cdot) \|_{L^5(\R^d)}  \ge \| \phi\|_{L^5(\R^d)}  -\left\|u(t,
      \cdot)-\phi\right\|_{L^5(\R^d)}\ge m_5 - \left\|u(t,
      \cdot)-\phi\right\|_{L^5(\R^d)},
  \end{equation*}
this implies that 
\[
\inf_{t \in \R} \| u(t, \cdot) \|_{L^5(\R^d)} \ge m_5-\eps >0.
\]
In particular, since $\|u(t, \cdot)\|_{L^2}={\rm const.}$, an
interpolation between $L^2$ and $L^\infty$ proves
\[
\inf_{t \in \R} \| u(t, \cdot) \|_{L^\infty(\R^d)} >0,
\]
and hence \eqref{eq:disp-Linfty} cannot hold.
\end{proof}

\subsection{Further mathematical results} 
In the following we review some of the known results on orbital (in-)stability of ground states in $d=1,2,3$. Moreover, we shall 
prove that 3D action  ground states are not
necessarily energy minimizers.


\subsubsection{The purely cubic case}
\label{sec:cubic}
For the cubic Schrödinger equation
\begin{equation}
  \label{eq:cubicNLS}
  i\partial_t u +\frac{1}{2}\Delta u = -|u|^2 u,\quad (t,x)\in \R\times\R^d,
\end{equation}
with $d\le 3$, the Cauchy problem is globally well-posed for $d=1$ (both in
$L^2(\R)$ and 
$H^1(\R)$, since the nonlinearity is $L^2$-subcritical), while finite
time blow-up is possible if $d=2$ or $3$, see e.g. \cite{CazCourant}. Regarding the solitary
waves, the analogue of \eqref{eq:soliton} is
\begin{equation}
  \label{eq:soliton-cubic}
  -\frac{1}{2}\Delta \phi +\omega \phi -|\phi|^2 \phi =0,\quad \phi\in
  H^1(\R^d)\setminus\{0\}. 
\end{equation}
The corresponding Pohozaev identities imply that such a non-trivial solution exists only if
$\omega>0$. Conversely, for any given $\omega>0$,
\eqref{eq:soliton-cubic} has a {\it unique} positive radial solution. As a
matter of fact, since the nonlinearity is homogeneous, the role of
$\omega$ is explicit: consider $ Q_{\rm cubic}$ the unique positive
radial solution in the case $\omega=1$,
\begin{equation}
  \label{eq:Qcubic}
  -\frac{1}{2}\Delta Q_{\rm cubic}+Q_{\rm cubic} -Q_{\rm
    cubic}^3=0,\quad x\in \R^d. 
\end{equation}
Then for any $\omega>0$,
\begin{equation*}
  \phi_\omega(x) := \sqrt \omega \, Q_{\rm cubic}(x\sqrt\omega)
\end{equation*}
is a positive radial solution to \eqref{eq:soliton-cubic}. By
uniqueness of such solutions, $\phi_\omega$ is an action  ground state,
minimizing
\begin{equation*}
  E_{\rm cubic}(\phi)+\omega M (\phi)= \frac{1}{2}\|\nabla \phi\|_{L^2(\R^d)}^2
  -\frac{1}{2}\|\phi\|_{L^4(\R^d)}^4 +\omega\|\phi\|_{L^2(\R^d)}^2.
\end{equation*}
In particular, we readily compute
\begin{equation*}
  \|\phi_\omega\|_{L^2(\R^d)} = \omega^{1/2-d/4}\|Q_{\rm cubic}\|_{L^2(\R^d)} .
\end{equation*}
Recalling the Grillakis-Shatah-Strauss stability criterion, this directly implies that
cubic ground states are orbitally stable in $d=1$, and unstable in
$d=3$ (in fact, we have strong instability by blow-up, see e.g. \cite{CazCourant}). 
The 2D case is $L^2$-critical and instability stems from the
fact that the cubic ground state has exactly {\it zero energy} (as seen from
\cite{Weinstein83}). Arbitrarily small perturbations can therefore make the
energy negative, which consequently leads to finite time blow-up of the associated solution $u$ by a standard virial
argument.  


\subsubsection{Cubic-quintic case in 2D}
\label{sec:rappel2D}
For $d=2$, it follows from the analysis in \cite{CaSp} that any
solitary wave has a mass larger than that of the cubic 
ground state, i.e. for any $\omega\in (0, \tfrac{3}{16})$, and any solution to
\eqref{eq:soliton}  
\begin{equation*}
  \|\phi\|_{L^2(\R^2)}>\|Q_{\rm cubic}\|_{L^2(\R^2)},
\end{equation*}
where $Q_{\rm cubic}$ is the radial, positive solution to
\eqref{eq:Qcubic}.
\smallbreak

Recall that $Q_\omega$ denotes the action ground
  state. Having in mind the Grillakis-Shatah-Strauss theory, the following asymptotic
results have been proved in \cite{CaSp,LewinRotaNodari20}:
for $\omega\approx 0$ or $\omega\approx \tfrac{3}{16}$, the map $\omega\mapsto M(Q_\omega)$ is
{\it increasing}. This implies orbital stability in the sense of
Definition~\ref{def:stability}, at least for some range of
the frequency $\omega$ close to the critical values. 
The numerical plots of $M(Q_\omega)$ given in \cite{LewinRotaNodari20} (see also Section \ref{sec:numground} below) suggest that
$\omega\mapsto M(Q_\omega)$ is
indeed increasing on the whole range $\omega\in (0,\tfrac{3}{16})$, and hence:
\begin{conjecture}\label{conj:2d}
In $d=2$, all cubic-quintic action  ground state
solutions are orbitally stable.  
\end{conjecture}

In Section \ref{sec:stab2D} we shall give further numerical evidence for this conjecture to be true, by 
performing several simulations of the time-evolution of perturbed action  ground states in 2D.


\subsubsection{Cubic-quintic case in 3D}
\label{sec:rappel3D}

As established in \cite{KOPV17,LewinRotaNodari20}, when $d=3$, it holds:
\begin{itemize}
\item[(i)] On the one hand, as $\omega\to 0$, it holds:
  \begin{equation*}
    M(Q_\omega) = \frac{1}{\sqrt\omega}M(Q_{\rm
      cubic})+\frac{\sqrt\omega}{2} \|Q_{\rm
      cubic}\|_{L^6(\R^3)}^6+\mathcal O\(\omega^{3/2}\),
  \end{equation*}
where $Q_{\rm cubic}$ is the positive, radial solution to \eqref{eq:Qcubic}.  
\item[(ii)] On the other hand:
  \begin{equation*}
    \lim_{\omega \to {3}/{16}} M(Q_\omega) = \lim_{\omega \to {3}/{16}}  \frac{\partial  M(Q_\omega)}{\partial \omega}=+ \infty.
  \end{equation*}
\end{itemize}
According to the Grillakis-Shatah-Strauss theory, this implies that
cubic-quintic action  ground
states in 3D are {\it unstable} near 
$\omega_{\rm min}=0$, and 
orbitally {\it stable} near $\omega_{\rm max}=\tfrac{3}{16}$. Numerical
plots in \cite{LewinRotaNodari20} show a U-shaped curve for
$\omega\mapsto M(Q_\omega)$. This suggests the existence of a unique
unstable branch and a unique stable branch. We shall numerically investigate the nature of instability in this case
in Section \ref{sec:stab3D}.
\smallbreak

Recalling the fact that the set of (constrained) energy minimizers with negative energy is indeed orbitally stable, cf. Theorem \ref{theo:stability-constrained}, we shall 
now show that solutions to \eqref{eq:soliton} may
have positive energy when $d=3$. Indeed, following the approach of
\cite[Section~2.4]{KOPV17}, we can rescale $Q_\omega$ via
\begin{equation*}
  \psi_\omega(x)
  :=\frac{1}{\sqrt\omega}Q_\omega\(\frac{x}{\sqrt\omega}\). 
\end{equation*}
The new unknown $\psi_\omega$ then solves
\begin{equation*}
  -\frac{1}{2}\Delta \psi_\omega +\psi_\omega -\psi_\omega^3+\omega
  \psi_\omega^5=0,
\end{equation*}
and, as established in \cite{KOPV17},
\begin{equation*}
  \psi_\omega =Q_{\rm cubic}+\mathcal O(\omega)\text{ in }H^1(\R^3),\quad
  \text{as }\omega\to 0. 
\end{equation*}
This implies
\begin{align*}
  E(Q_\omega) &= \sqrt\omega\(\frac{1}{2}\|\nabla Q_{\rm
    cubic}\|_{L^2(\R^3)}^2 - \frac{1}{2}\| Q_{\rm
                cubic}\|_{L^4(\R^3)}^4+\mathcal O(\omega)\) \\
  &= \sqrt\omega\(
  \|Q_{\rm cubic}\|_{L^2(\R^3)}^2 +\mathcal O(\omega)\), 
\end{align*}
where the last simplification stems from Pohozaev identities for
$Q_{\rm cubic}$ (discard the $L^6$ norms from \eqref{eq:phi1} and
\eqref{eq:phi2}). Therefore, there exists $\omega_0>0$ such that
\begin{equation*}
  E(Q_\omega)>0,\quad \forall \omega\in (0,\omega_0). 
\end{equation*}
Recalling that $M(Q_\omega)\to \infty$ as $\omega\to 0$, this shows
that there exists unstable action  ground states with positive energy and
arbitrarily large mass. On the other hand, Theorem~\ref{theo:stability-constrained}  shows that
there exists $\rho_0>0$ such that for all $\rho>\rho_0$, 
the minimization problem \eqref{eq:8.3.5} has a solution, and $E_{\rm min}(\rho)<0$. In summary this yields:

\begin{theorem}\label{theo:different}
  Not all action  ground states in $d=3$ are
  energy ground states. 
\end{theorem}

To our knowledge, this is the first rigorous statement
  which shows that the  
two notions of action versus energy ground states, in general, need to
be considered as independent (this is also a consequence of
\cite[Appendix~E]{Weinstein85}, as pointed out in
\cite{LewinRotaNodari20}, since at least for some values of $\omega$
close to zero, $\partial M(Q_\omega)/\partial \omega<0$).
Note that our theorem is consistent with the results of  \cite{JeanjeanLu-p} which establishes that the converse is true, i.e. all energy 
ground states are action ground states. From the point of view of dynamics, the results above leave open the possibility 
of having orbitally stable least action ground states, which are \emph{not} members of $\mathcal E(\rho)$.


\subsubsection{Previous results in the 1D case}

In the case $d=1$, it is fairly natural to generalize the nonlinearity in \eqref{eq:nls} to recover features similar to
those of \eqref{eq:nls} when $d=2$ or $3$. More precisely, consider
\begin{equation}
  \label{eq:1d}
  i\partial_t u +\frac{1}{2}\partial_x^2 u = -|u|^{p-1}u+|u|^{q-1}u,
\end{equation}
with $1<p<q$. In dimension one, all algebraic nonlinearities are
energy-subcritical and explicit solution formulas for $\phi_\omega$
are available in some cases (see in particular
  \cite{IlievKirchev93}).  
In addition, the focusing part is $L^2$-critical for $p=4$. We therefore expect \eqref{eq:1d} to 
behave similar to \eqref{eq:nls} in $d=2$, if we choose $p=4$, and
similar to \eqref{eq:nls} in $d=3$, if $p>4$. Indeed,
using \cite{GSS87,IlievKirchev93}, it is proved in
\cite{Ohta95} that for $p=4$, all standing waves (not
  necessarily ground states) are orbitally stable,
while for $p>4$, some are orbitally stable (for $\omega \approx
\omega_{\rm max}$, computed analogously to the value $\tfrac{3}{16}$
for \eqref{eq:nls}),
and some are unstable (for $\omega \approx 0$).

Numerical simulations have addressed the case $p>4$, see
\cite{BuGr01,Grikurov95,Sulem}. In particular, \cite{Grikurov95}
reports simulations for perturbations of unstable action
  ground states, showing two 
possible dynamics: full dispersion, or convergence to another (stable)
soliton. 

\begin{remark}
  In the case $d=1$ and $p>4$, the conclusion of
  Theorem~\ref{theo:different} remains true, using the same proof.
\end{remark}


\section{Numerical construction of action  ground states}\label{sec:numground}

\subsection{Numerical algorithm} In this section, we shall discuss a numerical approach for constructing least action ground state solutions to \eqref{eq:soliton} in dimensions $d=2$ and $3$.
To this end, we first note that since $Q_\omega$ is real and radially symmetric, it solves
\begin{equation} \label{Qeqr}
	\frac{1}{2}\left(\frac{\partial^2Q_\omega}{\partial r^2} +\frac{d-1}{r}\frac{\partial Q_\omega}{\partial r} \right)-\omega Q_\omega+Q_\omega^{3}-Q_\omega^{5}=0,
\end{equation}
where $r=|x|$. 
In order to get an equation with regular coefficients (which consequently allows for a more efficient numerical approximation), we introduce the new 
independent variable 
\begin{equation}
	s=r^{2},
	\label{s}
\end{equation}
in which (\ref{Qeqr}) reads
\begin{equation}\label{Qeqs}
	2s \frac{\partial^2Q_\omega}{\partial s^2} +d \frac{\partial Q_\omega}{\partial s}-\omega Q_\omega+Q_\omega^{3}-Q_\omega^{5}=0.
\end{equation}
Since it is known that cubic-quintic ground states are exponentially decreasing (see, e.g., \cite{CaSp}), we choose an $s_{0}\gg1$ such that $Q_\omega(s_{0})$ vanishes within 
numerical precision (which is of the order of $10^{-16}$ here since we work in double precision). Below $s_{0}=10^3$, while in the next section we shall also consider examples 
with $s_{0}=10^{4}$. The numerical task is thus to find a 
non-trivial solution to (\ref{Qeqs}) for given $\omega \in (0, \tfrac{3}{16})$, such that $Q_\omega$ 
(numerically) satisfies the homogenous Dirichlet condition 
$Q_\omega(s_{0})=0$. 

The interval $[0,s_{0}]$ is then mapped via $s = \frac{s_{0}}{2}(1+\ell)$, $\ell\in[-1,1]$ to the 
interval $[-1,1]$. On the latter we introduce standard \emph{Chebyshev 
collocation points} $\ell_{n}=\cos(n\pi/N$), $n=0,\ldots,N$, $N\in 
\mathbb{N}$ to discretize 
the problem. For any given $\omega>0$ in the admissible range, the function $Q\equiv Q_\omega$ is consequently approximated via the {\it Lagrange 
interpolation polynomial} $P_{N}(\ell)$ of degree $N$, coinciding with $Q$ at the collocation 
points, 
\[
P_{N}(\ell_{n}) = Q(\ell_{n}),\quad n = 0,\ldots,N.
\] 
Similarly, the (radial) derivative of $Q$ is approximated via the derivative of the Lagrange 
polynomial, i.e.
\[
\frac{\partial}{\partial s} Q(s(\ell_n))\approx P_{N}'(\ell_n).
\] 
At the collocation points, this implies $\partial_sQ(\vec{\ell})\approx D \mathrm{Q}$,  
since the interpolation polynomial is obviously linear in the $\ell_{n}$, 
$n=0,\ldots,N$. Here, $\vec{\ell}$ is the vector with components $\ell_n$, $D$ is the \emph{Chebyshev 
differentiation matrix} \cite{trefethen,WR} and $\mathrm{Q}$ is the vector with 
components $Q(s(\ell_{n}))$, $n=0,\ldots,N$. 

With the above discretization, equation (\ref{Qeqs}) is approximated 
by a system of nonlinear equations for the vector $\mathrm{Q}$ which 
can be formally written in the form $\mathbb{F}(\mathrm{Q})=0$. The homogenous Dirichlet condition 
for $s=s_{0}$ is thereby implemented by 
eliminating the column and the line corresponding to $s_{0}$, cf. \cite{trefethen} for more details. This
shows that $\mathbb{F}(Q)=0$ is an $N$-dimensional system of 
nonlinear equations for the $N$ components $Q(s(\ell_{n}))$, 
$n=1,\ldots,N$. 
This system will be solved via a Newton iteration,
\begin{equation}
	\mathrm{Q}^{(m+1)} = \mathrm{Q}^{(m)}- 
	(\mathrm{Jac}\, \mathbb{F}(\mathrm{Q}^{(m)}))^{-1}\mathbb{M}\mathrm{Q}^{(m)}
	\label{newton},
\end{equation}
where $\mathrm{Jac}\, \mathbb{F}$ denotes the Jacobian of
$\mathbb{F}$, and where  
$\mathrm{Q}^{(m)}$, $m=0,1,\ldots$ denotes the $m$-th iterate. 

Note, however, that $\mathrm{Q}=0$ is always a trivial solution to
\eqref{Qeqs}, which needs to be avoided during the iteration
process.  
In order for this algorithm to converge to our desired, non-trivial solution $Q$, one needs to identify a suitable initial iterate $\mathrm{Q}^{(0)}$. To do so, we shall apply a 
tracing or continuation technique as follows: We introduce in (\ref{Qeqs}) a deformation parameter 
$\alpha\in [0,1]$, such that for $\alpha=0$ we have only the focusing cubic 
nonlinearity, while for $\alpha=1$ we obtain the full cubic-quintic equation, i.e. we effectively solve 
\begin{equation}\label{Qeqsa}
	2s \frac{\partial^2Q_{\omega,\alpha}}{\partial s^2} +d \frac{\partial Q_{\omega,\alpha}}{\partial s}-\omega Q_{\omega,\alpha}+Q_{\omega,\alpha}^{3}-\alpha Q_{\omega,\alpha}^{5}=0,\quad 
	\alpha\in[0,1],
\end{equation}
instead of only \eqref{Qeqs}.
The cubic solitons $Q_{\omega,\alpha=0}$ are numerically known, see, e.g. \cite{AKS, KMS} (no explicit ground state formula exists in dimensions $d\ge 2$). 
We can thus solve the discretized 
equation \eqref{Qeqsa} for $\alpha=0.1$ and, say, $\omega=0.1$ via the Newton iteration described above. The 
resulting solution $Q_{\omega=0.1, \alpha=0.1}$ is then used as an initial iterate for the same 
equation for $\alpha=0.2$, and so on, until the cubic-quintic case $\alpha=1$ is reached. 
In a second step, we use the ground states obtained for $\omega=0.1$ as an initial iterate 
for slightly smaller or larger $\omega$'s within the admissible range $0<\omega < \frac{3}{16} =0.1875$. In this way all examples presented 
below can be treated. 

\begin{figure}[htb!]
  \includegraphics[width=0.49\textwidth]{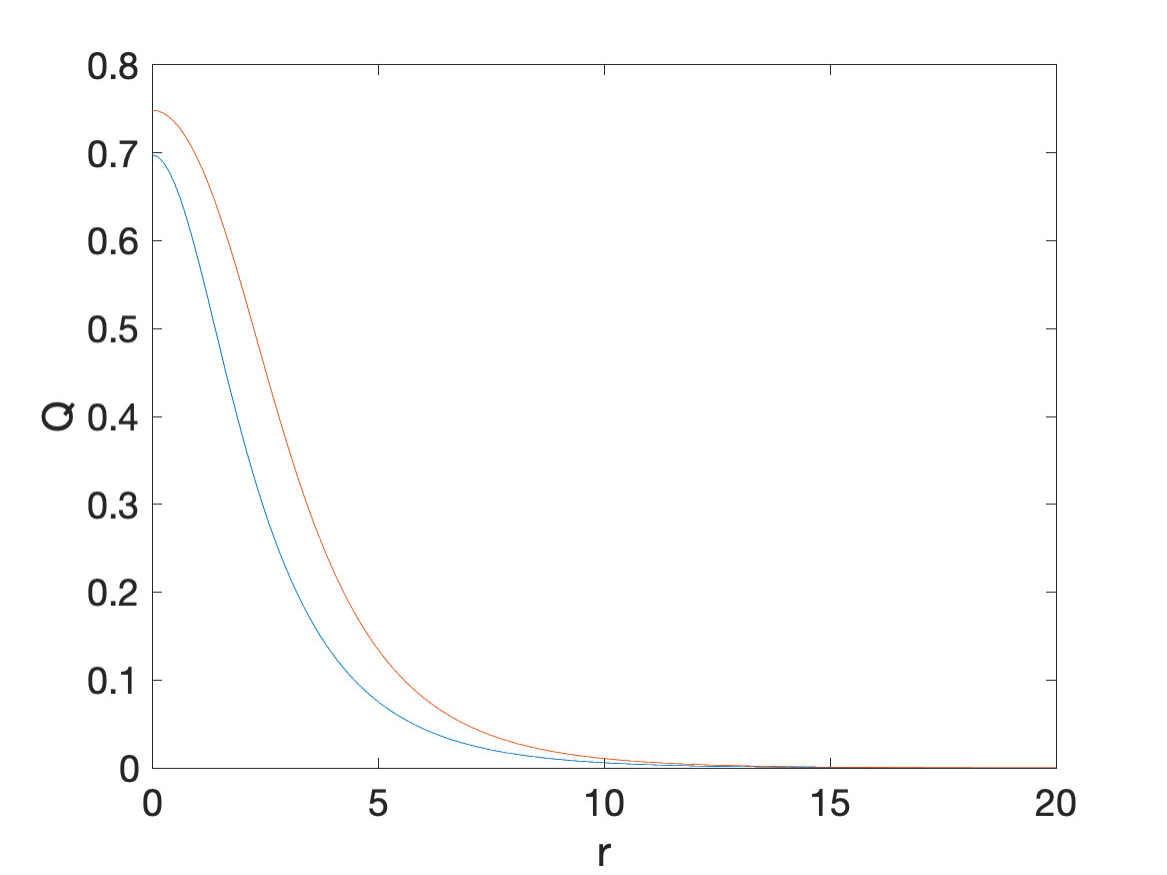}
  \includegraphics[width=0.49\textwidth]{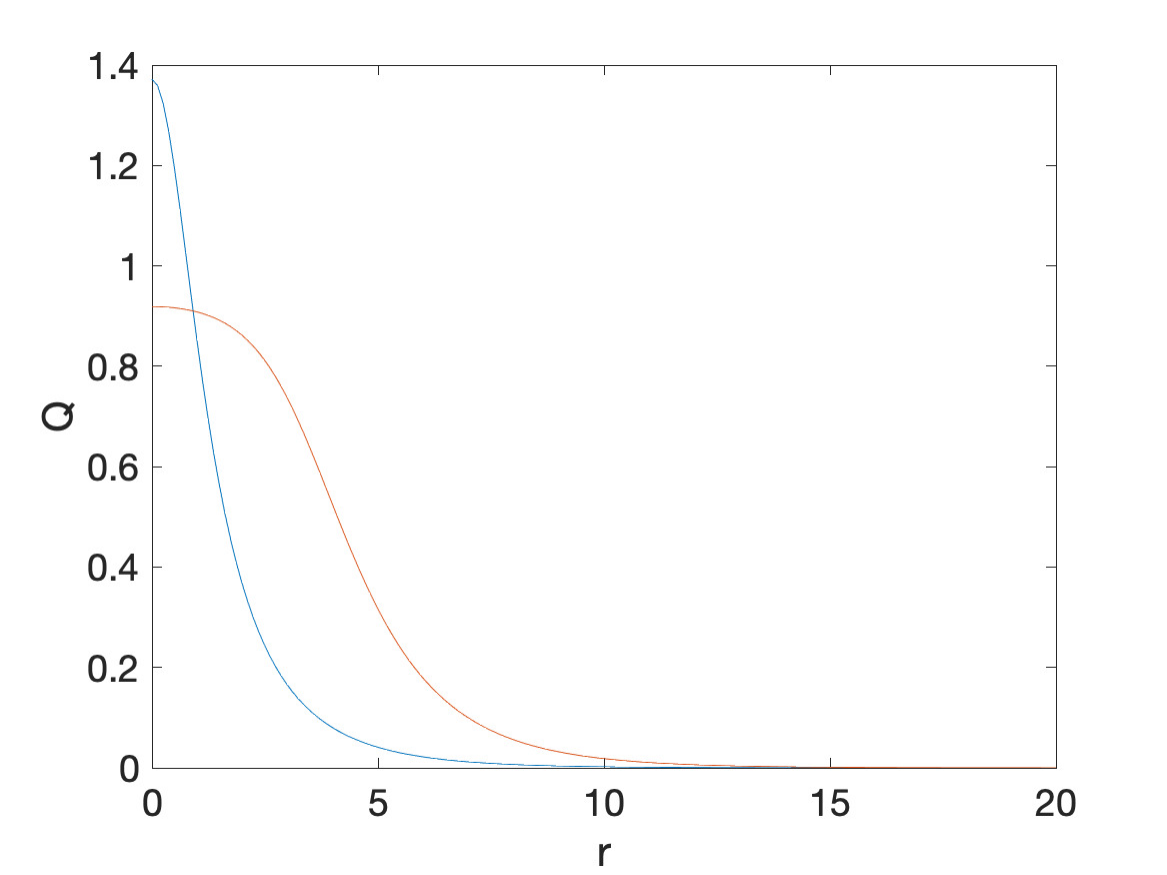}
\caption{Ground state solutions $Q_{\omega=0.1}$ to the cubic NLS in 
blue and the cubic-quintic NLS in red: on the left for $d=2$ and on the 
right for $d=3$.}
 \label{2d3d}
\end{figure}

As an example we show in Fig.~\ref{2d3d} the ground state solutions $Q_\omega (r)$ 
at $\omega=0.1$
for the cubic NLS in blue and the cubic-quintic NLS in red.  
It can be seen that the situation is qualitatively different depending on the 
spatial dimension. Whereas in 2D, the cubic-quintic ground state has a slightly 
greater maximum and is slightly faster decaying than $Q_{\rm cubic}$, 
in 3D the cubic-quintic ground state has a much 
smaller maximum and a considerably larger support.


\subsection{Numerical ground states in 2D}
\begin{figure}[htb!]
  \includegraphics[width=0.49\textwidth]{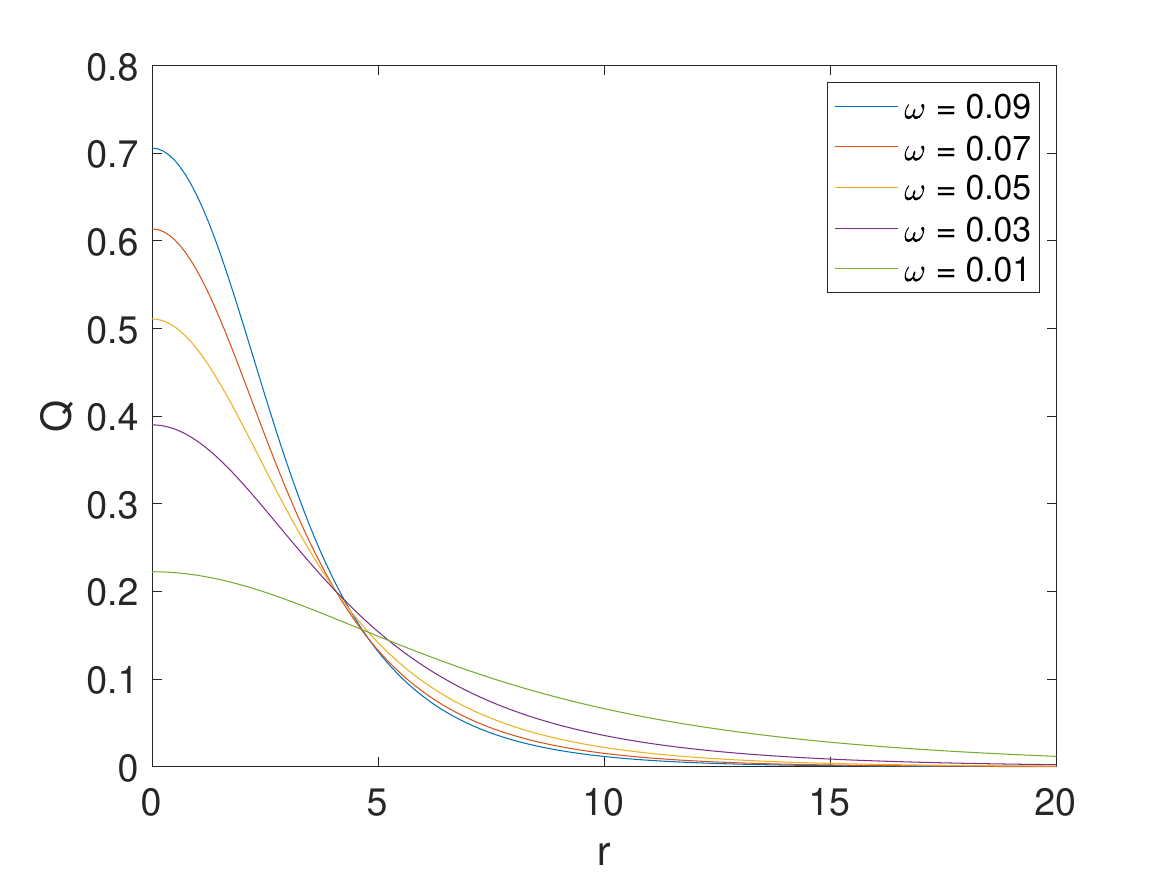}
  \includegraphics[width=0.49\textwidth]{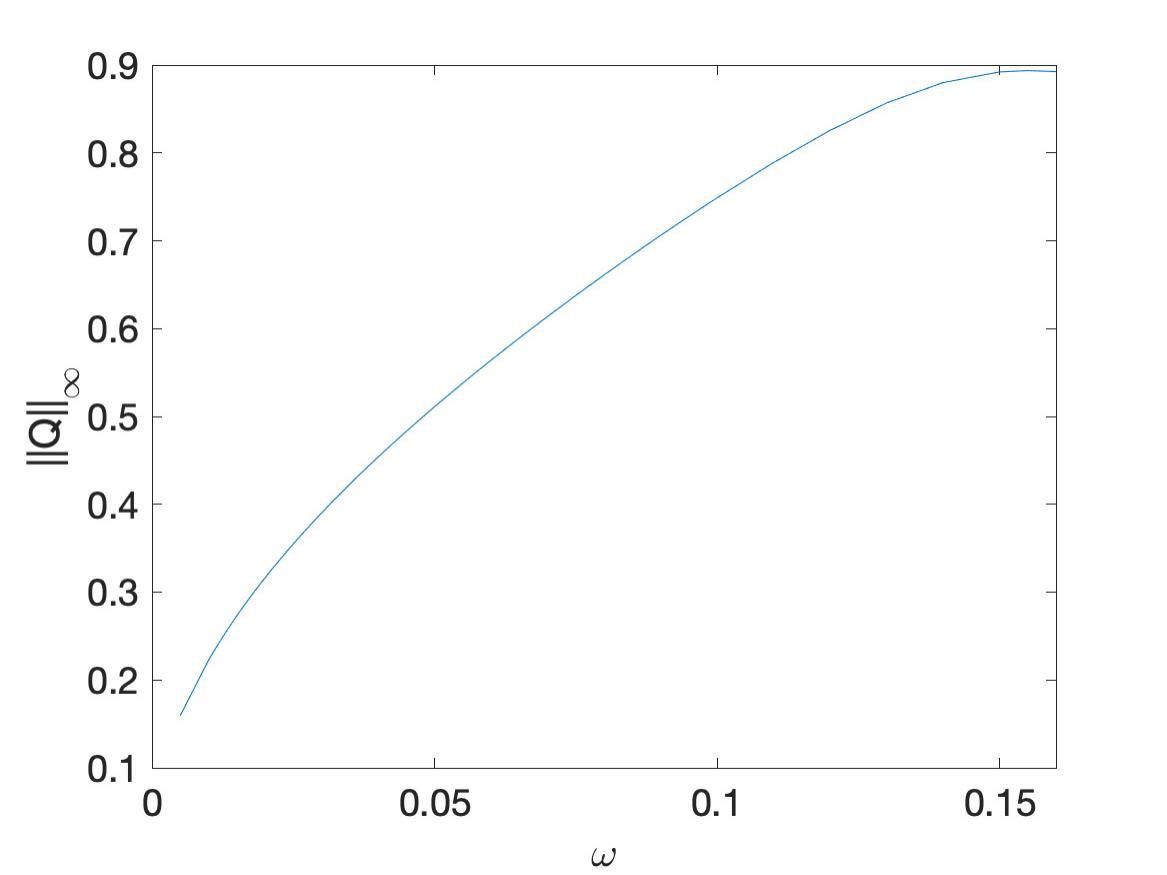}
\caption{Left: Ground state solutions to \eqref{Qeqr} in dimension $d=2$
 for several values of $\omega$. Right: The $L^{\infty}$-norm of these 
 states as a function of $\omega$.}
 \label{NL35_d2sol}
\end{figure}

We first consider the case $d=2$ with $N=400$ collocation points: In Fig.~\ref{NL35_d2sol} we show on 
the left a plot of the ground state function $Q_\omega(r)$ for various values of $\omega$. It is seen that the maximum of 
the ground states increases with $\omega$. The solutions also become 
more localized with increasing $\omega$. On the right of the same 
figure, we show the $L^\infty(\R^2)$-norm of the ground states as function of  
$\omega$. For convenience, we only consider values of 
$\omega\in[0.005,0.16]$. 

\begin{figure}[htb!]
  \includegraphics[width=0.49\textwidth]{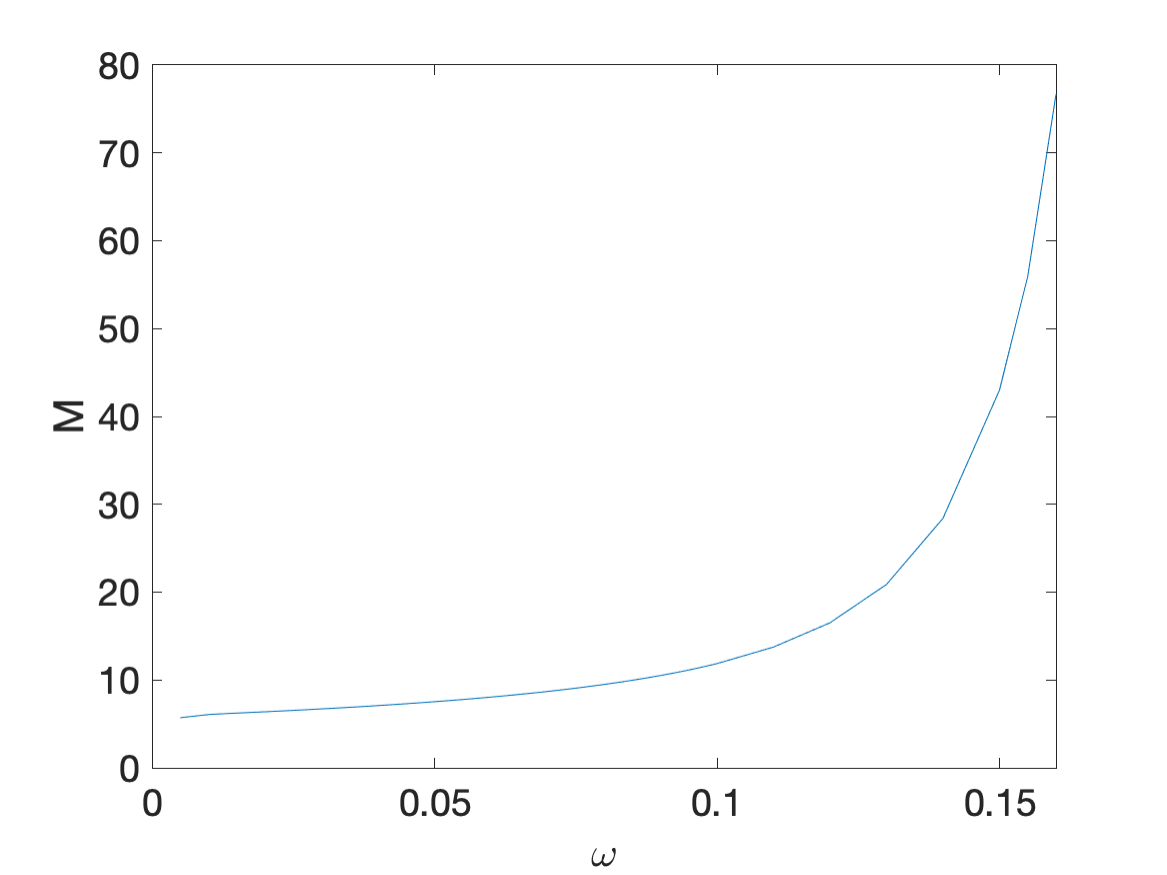}
  \includegraphics[width=0.49\textwidth]{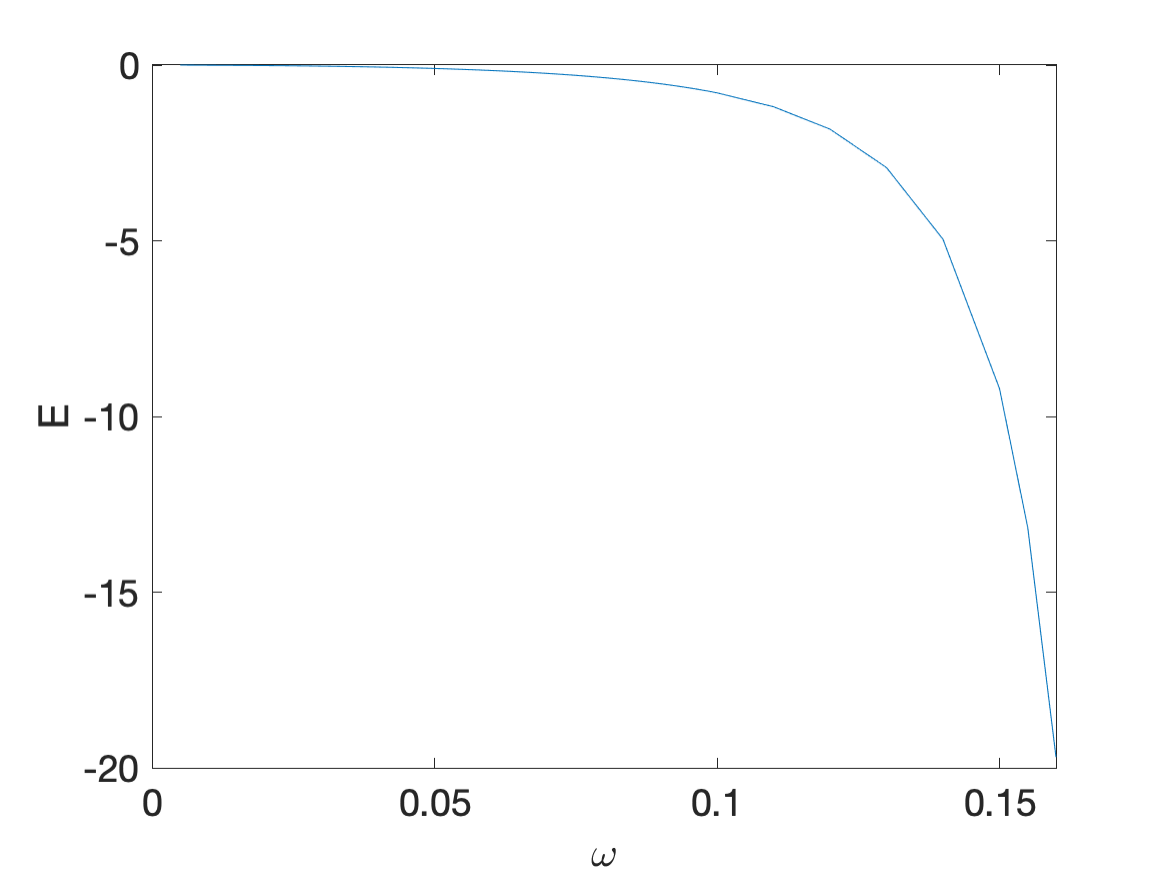}
 \caption{$M(Q_\omega)$ and $E(Q_\omega)$ as a functions of $\omega$ in dimension $d=2$.}
 \label{NL35_d2solmass}
\end{figure}

In Fig.~\ref{NL35_d2solmass} we depict 
the ground-state mass $M(Q_\omega)$ and energy $E(Q_\omega)$ as a
function of $\omega$. These plots are based on a total library 
of roughly 100 numerical ground state solutions $Q$ on the shown 
range of $\omega$.  The corresponding mass- and energy-integrals are 
thereby computed with the {\it Clenshaw-Curtis algorithm} in $s=r^2$, a spectral integration method 
based on the same Chebyshev collocation points as before, see \cite{trefethen}. 
Both $M(Q_\omega)$ and $E(Q_\omega)$ appear to be monotonic 
in $\omega$. In particular, the monotonicity of $M(Q)$ indicates
orbital stability in the sense of Definition~\ref{def:stability}, in
view of the Grillakis-Shatah-Strauss theory.

\subsection{Numerical ground states in 3D} 
\begin{figure}[htb!]
  \includegraphics[width=0.49\textwidth]{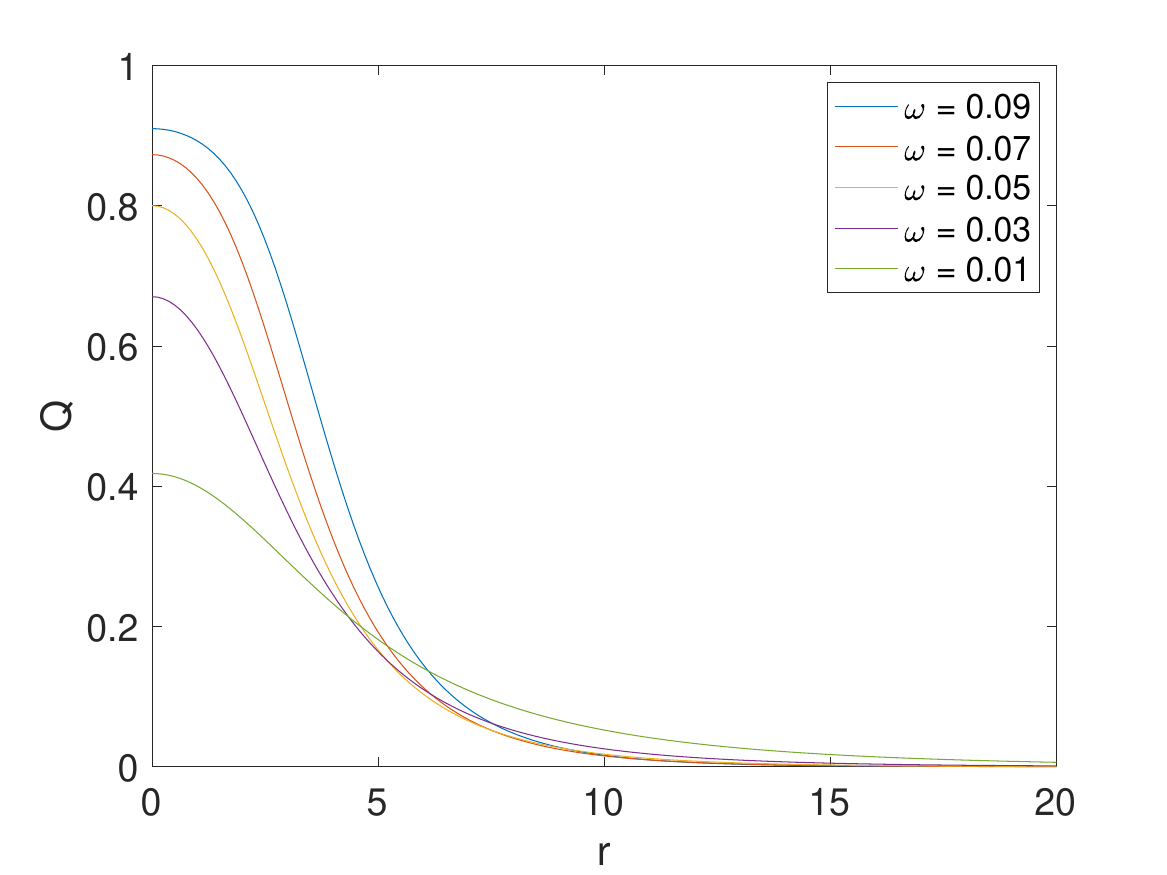}
  \includegraphics[width=0.49\textwidth]{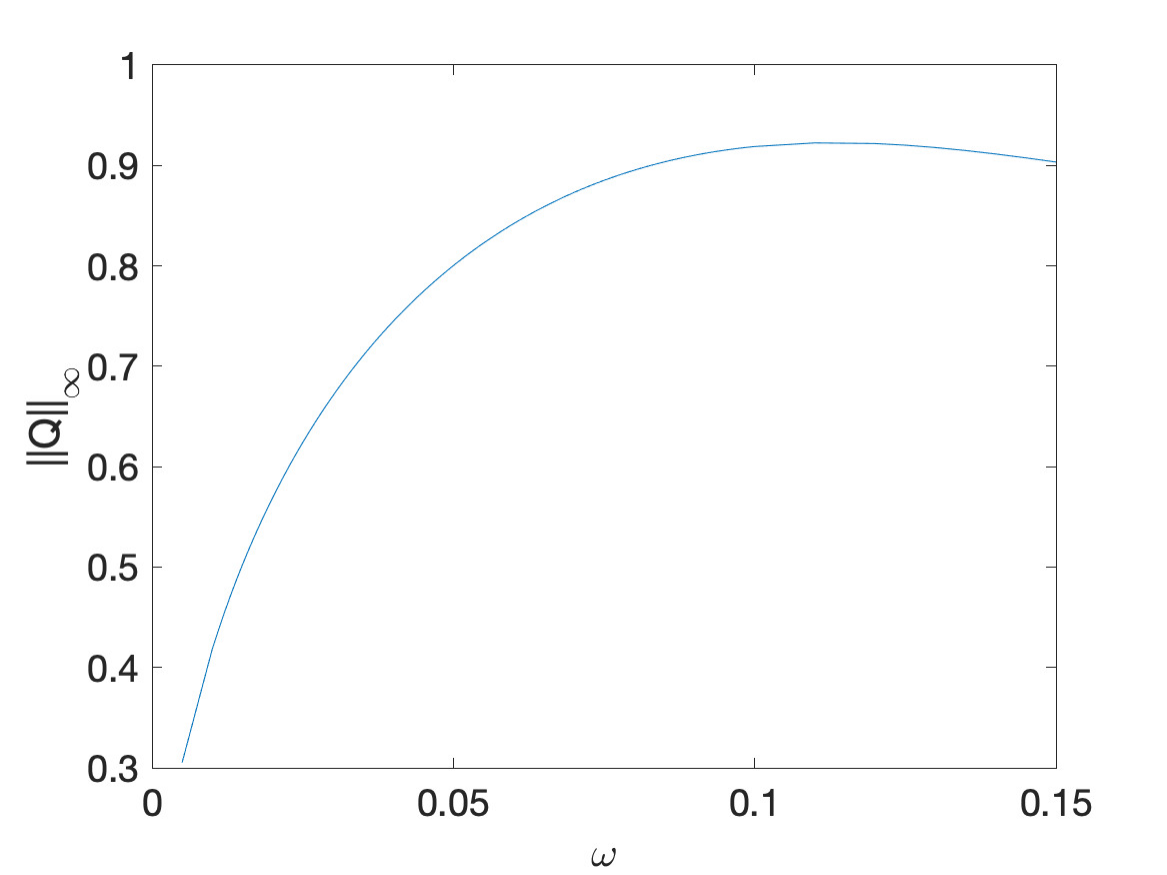}
 \caption{Left: Ground state solutions to \eqref{Qeqr} in dimension $d=3$
 for several values of $\omega$. Right: The $L^{\infty}$-norm of these 
 states as a function of $\omega$.}
 \label{NL35_d3sol}
\end{figure}

In the case $d=3$, we use the same numerical parameters as before: 
In Fig.~\ref{NL35_d3sol}, we show 
on the left the ground states for several values of $\omega$. It can 
again be seen that the maximum of $Q_\omega$ increases with $\omega$, 
at least up to some value $\omega_\ast\approx 0.1$. For larger values of $\omega$, however, the 
$L^\infty(\R^3)$-norm of $Q_\omega$ is seen to be decreasing again. 

Analogously to the 2D case, the solutions become more localized with increasing $\omega$. Note, 
however, that despite 
its exponential decay, the 3D soliton is less localized than in the 
case of the purely focusing, cubic NLS, see Fig.~\ref{2d3d} on the 
right. The 3D ground states in Fig.~\ref{NL35_d3sol} on the left
are also found to be less peaked 
than the corresponding solutions in dimension 2, see
Fig.~\ref{NL35_d2sol}. 

\begin{figure}[htb!]
 \includegraphics[width=0.49\textwidth]{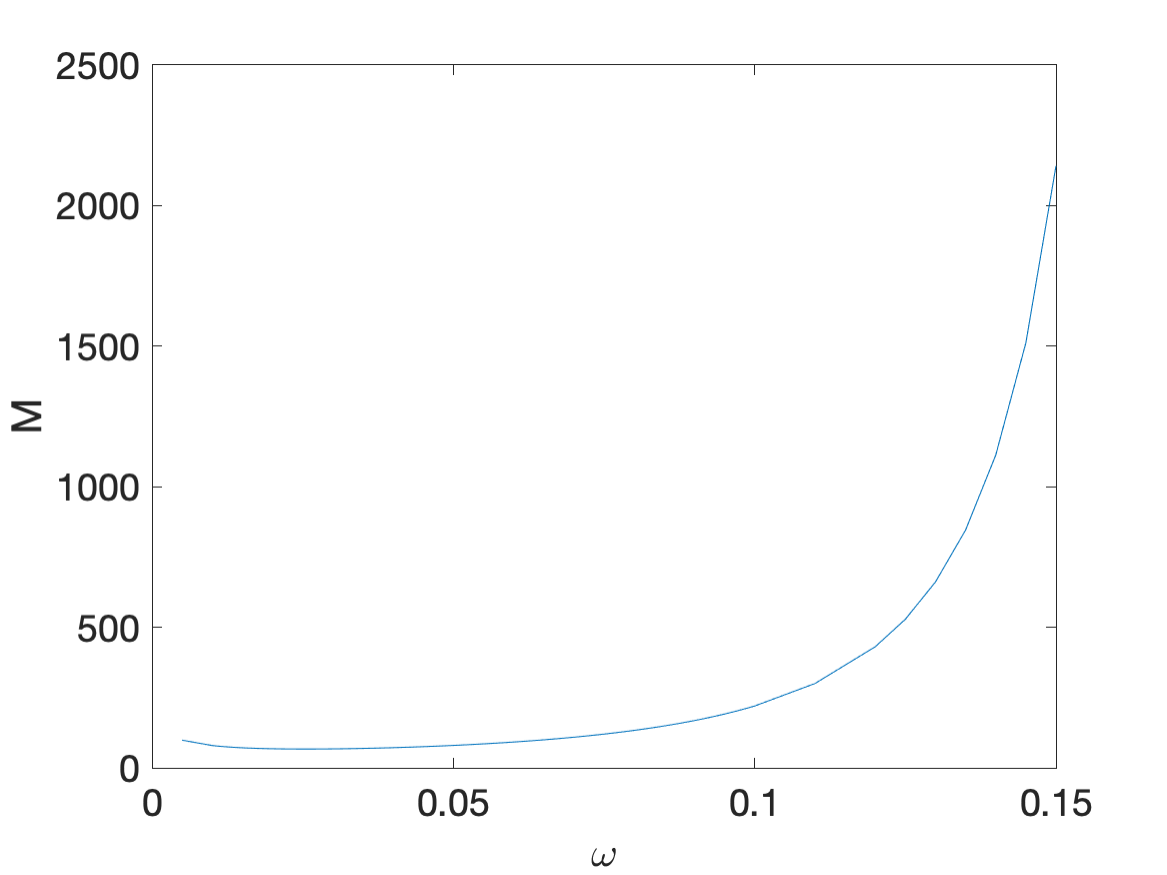}
 \includegraphics[width=0.49\textwidth]{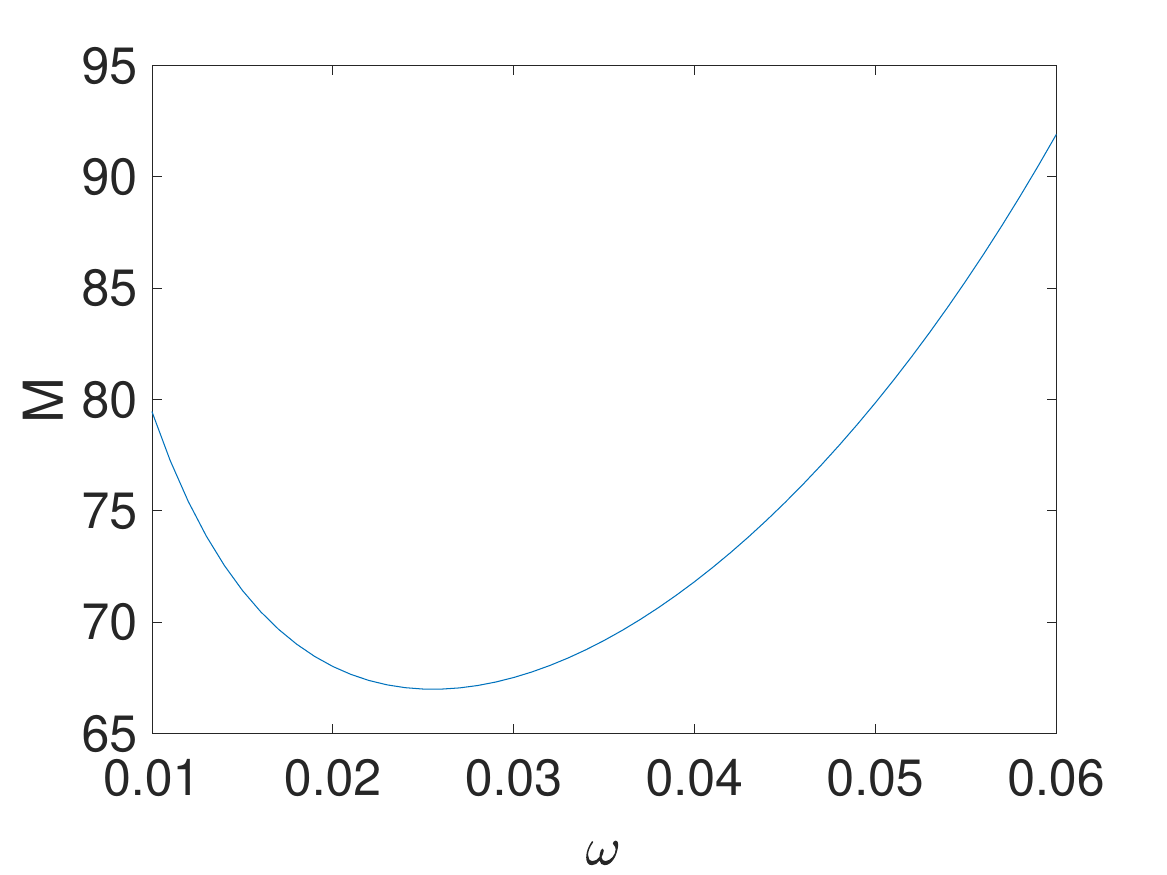}
 \caption{Left: $M(Q_\omega)$ as a functions of $\omega$, for cubic-quintic ground states $Q_\omega$ in dimension $d=3$. Right: a close-up 
 of the same curve near $\omega_{\rm crit}$.}
 \label{NL35_d3solmass}
\end{figure}

\begin{figure}[htb!]
  \includegraphics[width=0.49\textwidth]{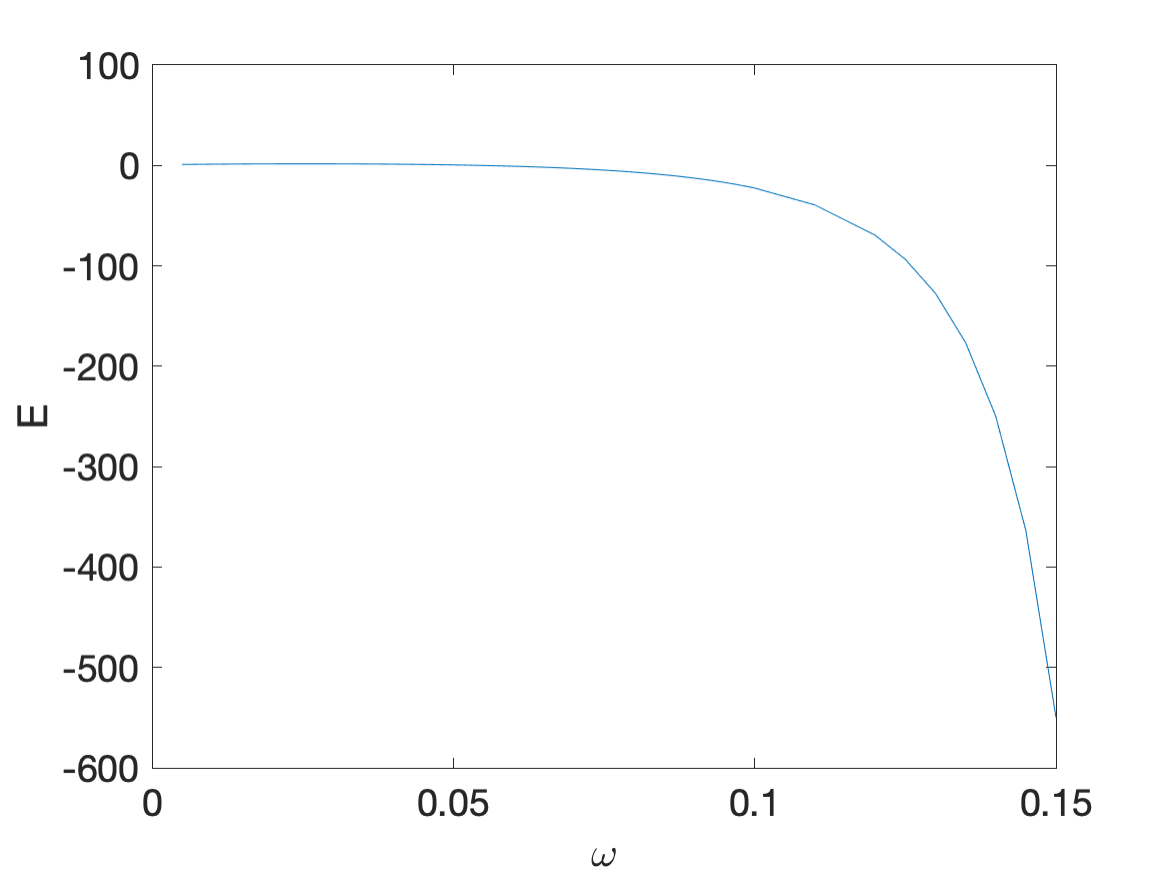}
  \includegraphics[width=0.49\textwidth]{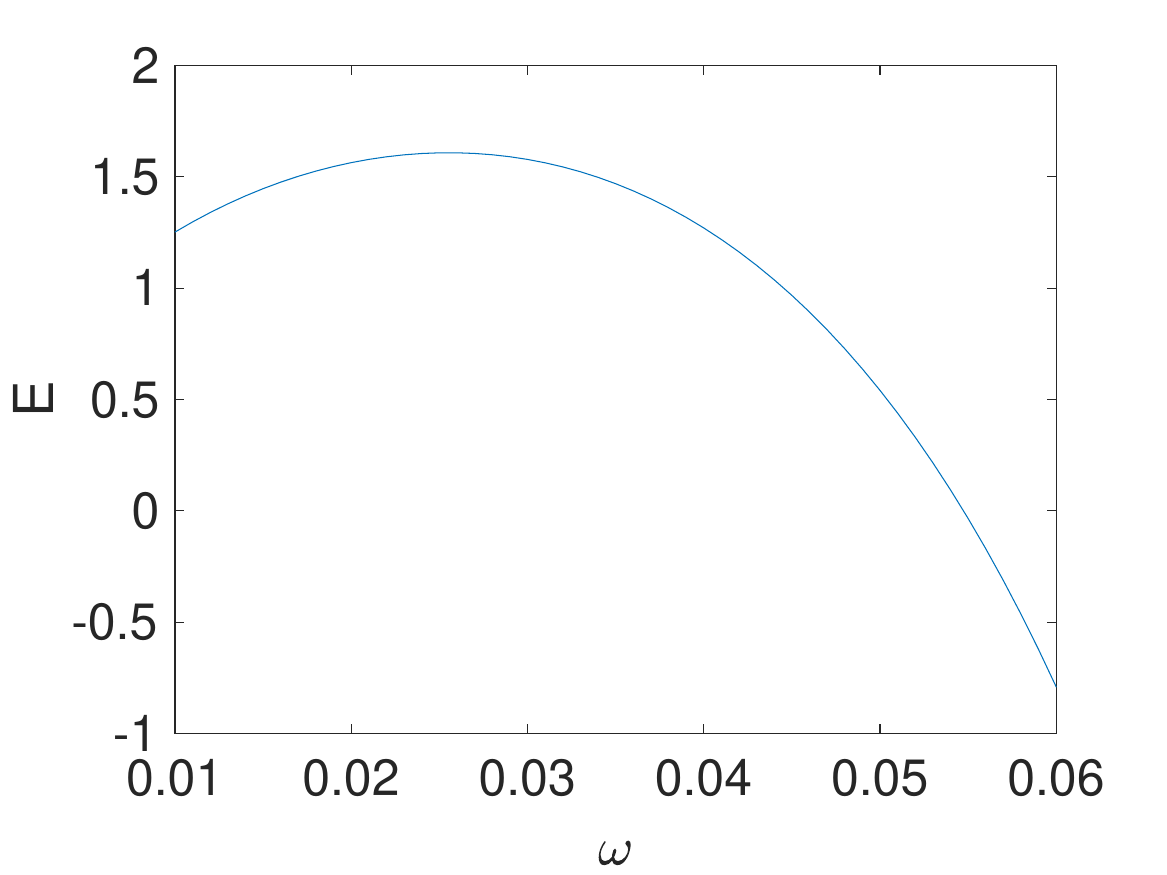}
 \caption{Left: $E(Q_\omega)$ as a functions of $\omega$, for cubic-quintic ground states $Q_\omega$ in dimension $d=3$. 
 Right: a close-up 
 of the same curve near $\omega_{\rm crit}$.}
 \label{NL35_d3solenergy}
\end{figure}

In contrast to the 2D case, the ground state mass $M(Q_\omega)$ is no longer monotonically increasing as a function of $\omega$.
Looking at Fig.~\ref{NL35_d3solmass}, we see that, instead, $M(Q_\omega)$ has a minimum at $\omega_{\rm crit}\approx 0.026$. 
We consequently expect orbital instability of ground states $Q_\omega$ for  
$\omega<\omega_{\rm crit}$, a phenomenon we shall study in more detail in Section \ref{sec:stab3D}.
In Fig. \ref{NL35_d3solenergy}, the corresponding ground state energy $E(Q_\omega)$ is seen to have a maximum 
at the same $\omega_{\rm crit}\approx 0.026$.

\begin{figure}[htb!]
  \includegraphics[width=0.49\textwidth]{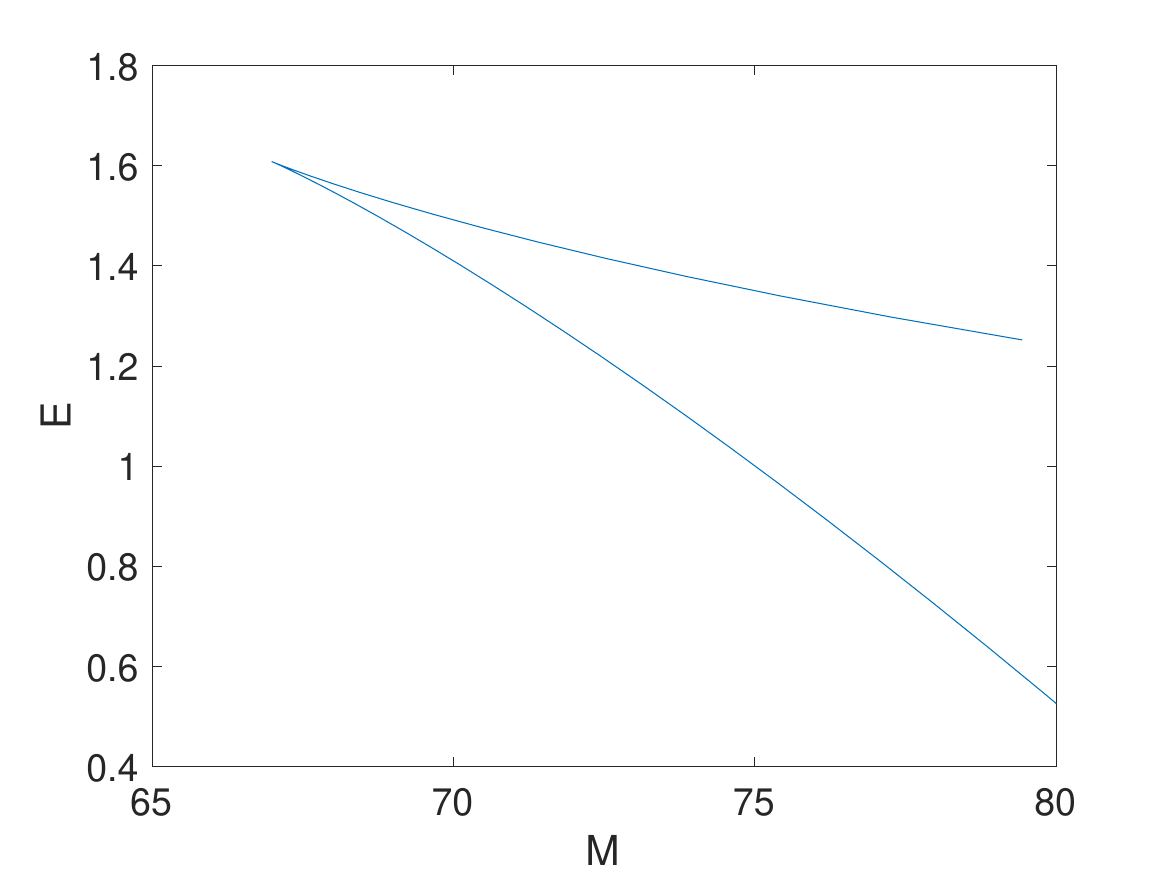}
 \caption{$E(Q_\omega)$ as a function of $M(Q_\omega)$ for cubic-quintic ground states in dimension $d=3$.}
 \label{NL35_d3Enermass}
\end{figure}

The appearance of an unstable branch is clearly visible when the energy
$E(Q_\omega)$ is plotted as a function of the mass $M(Q_\omega)$, see
Fig.~\ref{NL35_d3Enermass}. 
Our figure is in good agreement with
\cite[Figure~2]{KOPV17} (where the constants are different because
the factors are different from \eqref{eq:nls}). Clearly, 
ground states $\phi_\omega$ corresponding to the upper branch cannot correspond 
to constrained energy minimizers with mass $\rho = M(Q_\omega)$.


\section{Orbital stability of action  ground states in 2D}\label{sec:stab2D}

\subsection{Numerical method for the time-evolution}

In this section we will numerically study the time-evolution of \eqref{eq:nls} resulting from initial data $u_0$ given by 
perturbations of action  ground states. We will only consider 
perturbations which conserve the radial symmetry. To 
consider more general perturbations, a full 3D code would be 
necessary which is beyond the scope of the present paper. This allows us to use the change of variables \eqref{s} and effectively 
solve \eqref{eq:nls} in the (non-singular) form
\begin{equation}	\label{NLSpqs}
	i\partial_{t} u +2s \frac{\partial^2 u}{\partial s^2} +d \frac{\partial u}{\partial s} +|u|^{2}u-|u|^{4}u = 0, \quad d=2,3.
\end{equation}
We thereby use the same discretization for $s\in [0, \infty)$ in terms of Chebyshev 
collocation points as detailed in the previous section. 

After the spatial discretization in $s$, equation (\ref{NLSpqs}) is then approximated 
via a system of ordinary differential equations. These equations are then integrated in time using a
{\it time-splitting method} in which the linear step is solved numerically
via an implicit {\it fourth-order} Runge-Kutta method, see \cite{KS} for more details. The accuracy of this 
time-integration algorithm is henceforth controlled via the analytically conserved quantity $E(u)$, which 
in our case nevertheless depends on time due to unavoidable numerical errors. 
As discussed in \cite{etna}, the numerical conservation of the 
relative mass tends to overestimate the numerical error by one to 
two orders of magnitude. We shall always aim at a numerical error below 
$\mathcal O(10^{-3})$, i.e. below plotting accuracy. This means that we 
ensure a relative energy-conservation 
\[
\Delta_E= \left|\frac{E(t)}{E(0)}-1\right| 
\]
of order $\Delta_E= \mathcal O(10^{-5})$, or better.
 
We shall use a single computational domain $\Omega=[0, s_0]$ for which we impose a homogeneous 
Dirichlet condition $u(t,s_0)=0$, for all $t\ge 0$. We mostly choose $s_{0}=10^{3}$, 
but in some unstable situations we shall also take
$s_{0}=10^{4}$.  As a basic test case,
we first propagate the three-dimensional 
ground state $Q_\omega$, numerically found at $\omega=0.1$. We thereby use 
$N_t=10^{4}$ time-steps until a final time $t_{\rm f}=10$. We find that the hereby obtained numerical solution $u$, 
at $t=t_{\rm f}$, satisfies
\[
\max_{\Omega} | u (t_{\rm f}, \cdot) - e^{ i t_{\rm f} \omega}Q_\omega | =  \mathcal O(10^{-9}),
\] 
i.e. the same order of accuracy as reached in \cite{KS}.  

\begin{remark} In general, it is not 
unproblematic to work with a homogeneous Dirichlet boundary condition
on a finite  
numerical domain, since this could lead to unwanted reflections of the 
emitted radiation at the boundary, see, e.g., the discussion in 
\cite{birem}. In our case, however, only small, rapidly decreasing perturbations of 
ground states are considered. It is thus possible to work 
on sufficiently large computational domains $\Omega$, on which the radiation can separate from 
the bulk before spurious reflections from the boundary lead to 
noticeable effects. 
\end{remark}


\subsection{Time-evolution of perturbed 2D ground states} 

In this subsection, we shall study the time-evolution of perturbed ground states to \eqref{eq:nls} in dimension $d=2$. To this end, we first consider the case where
\begin{equation}
	u_0(x) = \lambda Q_\omega(x),\quad \lambda > 0.
	\label{ini}
\end{equation}
Here $\lambda>0$ is a perturbation parameter and $Q_\omega$ is a numerically obtained action  ground state, at a certain admissible frequency $\omega\in (0,\frac{3}{16})$. 

\begin{figure}[htb!]
  \includegraphics[width=0.49\textwidth]{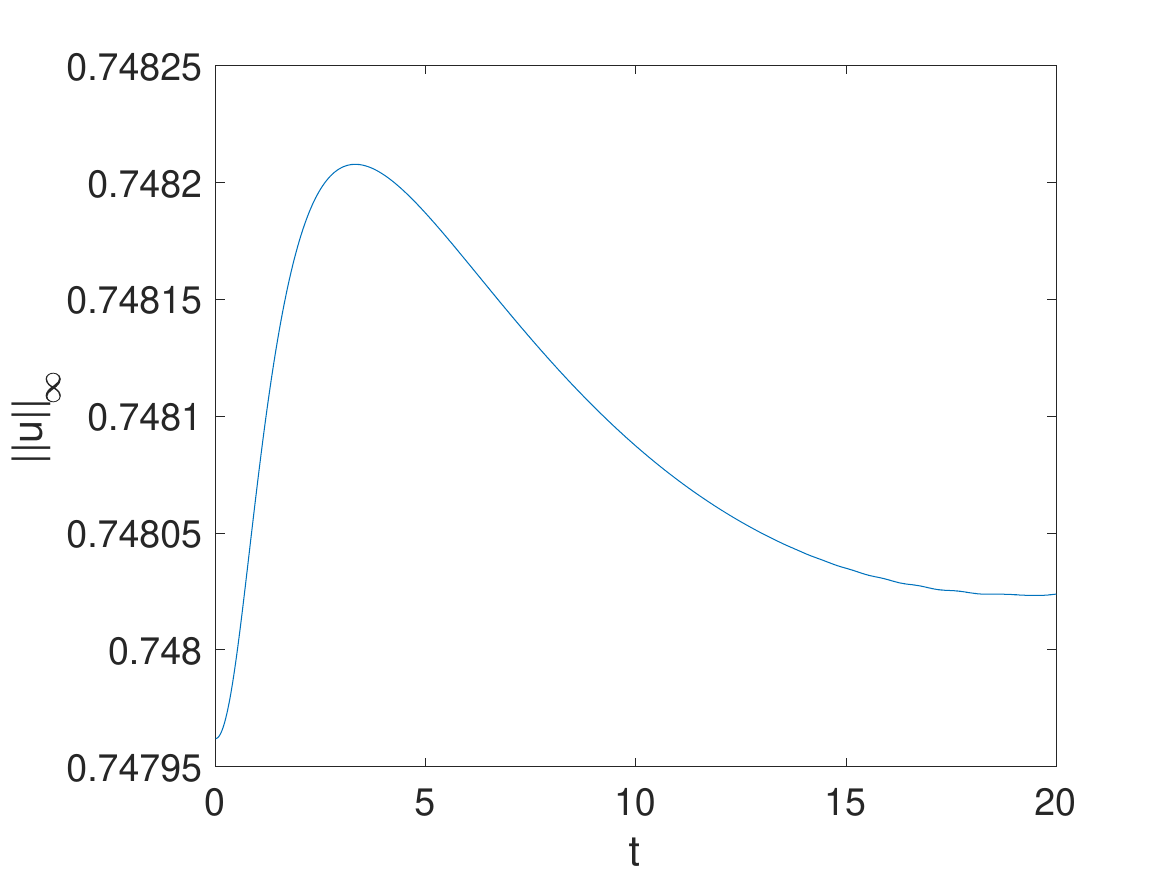}
  \includegraphics[width=0.49\textwidth]{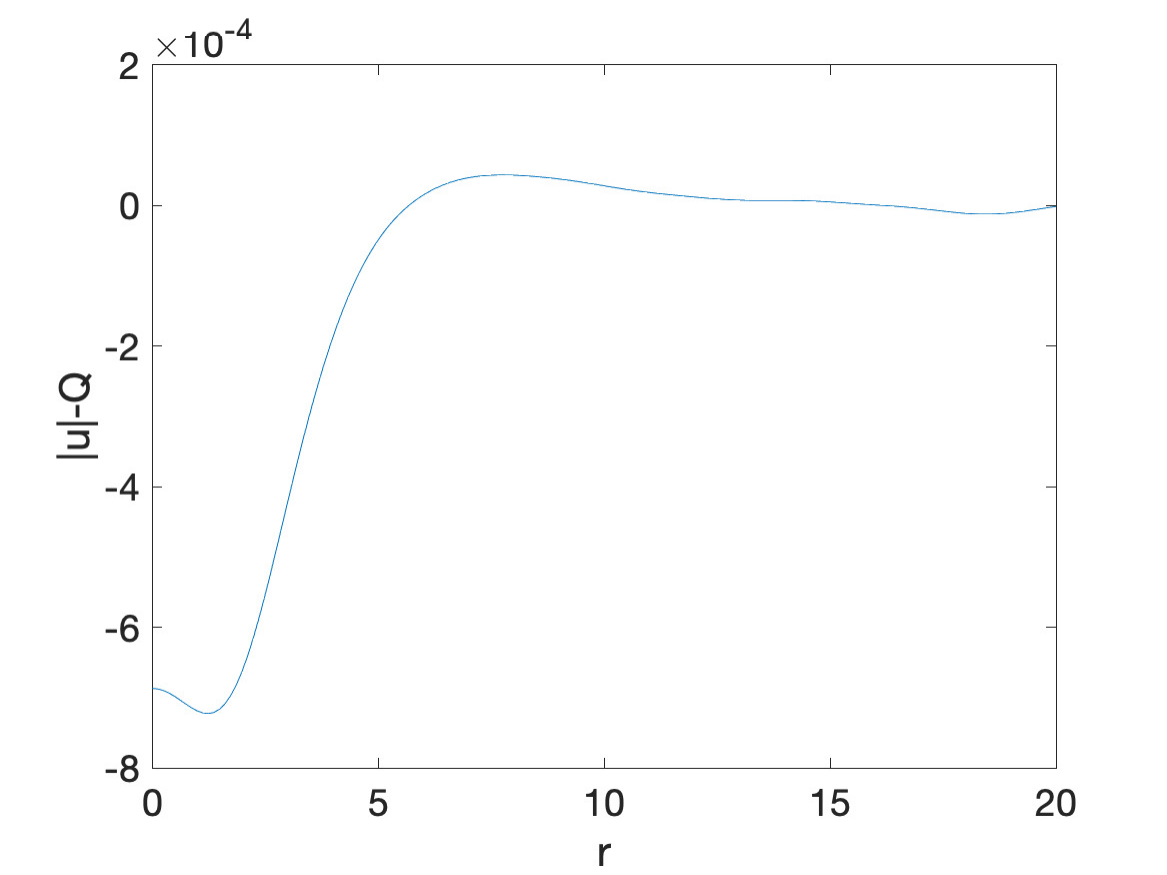}
 \caption{The 2D solution $u$ to \eqref{eq:nls} for initial data 
 (\ref{ini}) with $\omega=0.1$ and $\lambda =0.99$. On the left the $L^{\infty}$-norm as a function of time. On the right the difference between $|u|$ and $Q_{\rm \omega=0.1}$ at the final time $t_{\rm f}=20$.}
 \label{NL35_d2solom01}
\end{figure}

We first study the case where $\omega=0.1$ and $\lambda=0.99$, and use 
$N_{t}=10^{4}$ time steps to reach the indicated final time $t_{\rm f}=20$. As expected, the solution $u$ to \eqref{eq:nls}, effectively given by \eqref{NLSpqs}, 
is found to be close to the exact time-periodic
state \[\phi_\omega(t,x) = e^{i \omega t} Q_\omega(x).\] 
To this end, we show on the left of 
Fig.~\ref{NL35_d2solom01} the $L^{\infty}(\R^2)$-norm of the solution 
as a function of time. It can be seen that it approaches a final
state as $t\to 20$. The latter is found to be very close (in absolute value) to the unperturbed ground state $Q_{\rm \omega=0.1}$. 
Note that the $L^\infty$-difference is of the order of $10^{-4}$ and thus, 
much smaller than the initial perturbation. 

\begin{figure}[htb!]
  \includegraphics[width=0.49\textwidth]{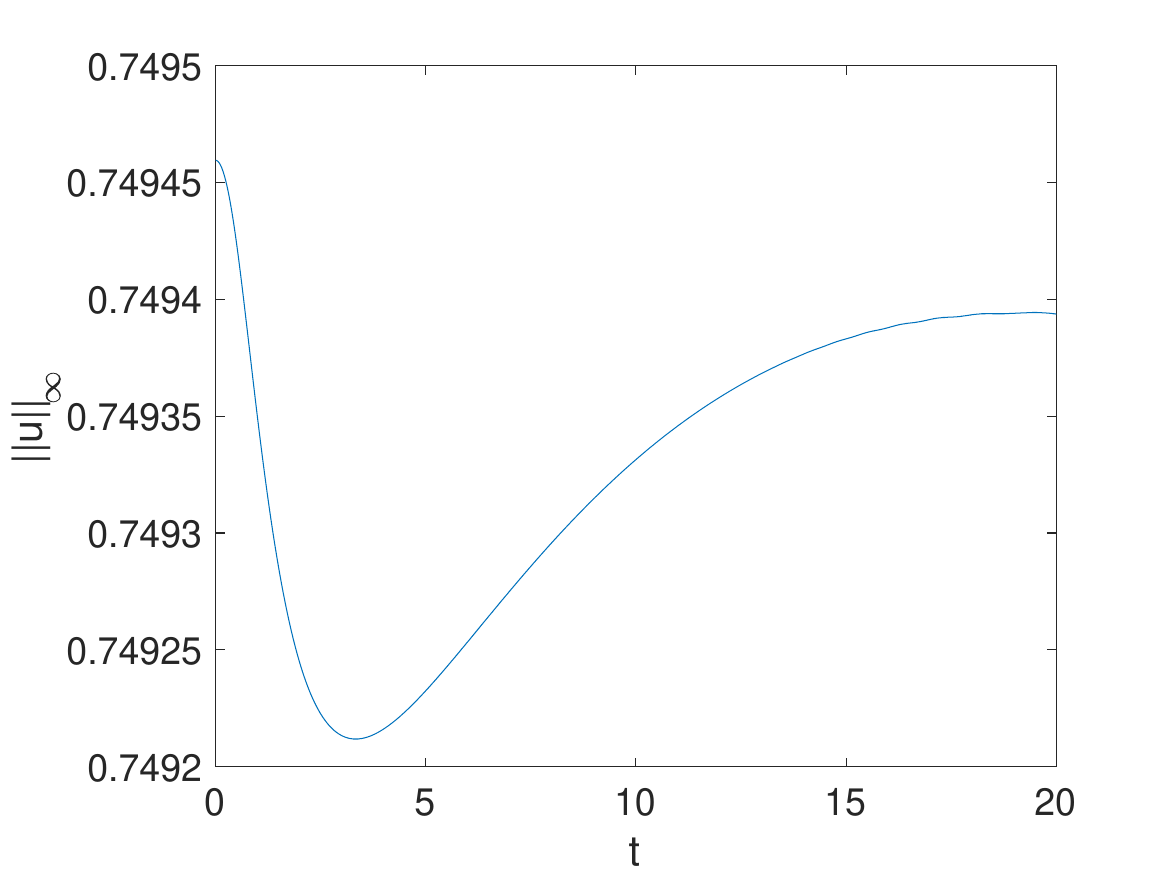}
  \includegraphics[width=0.49\textwidth]{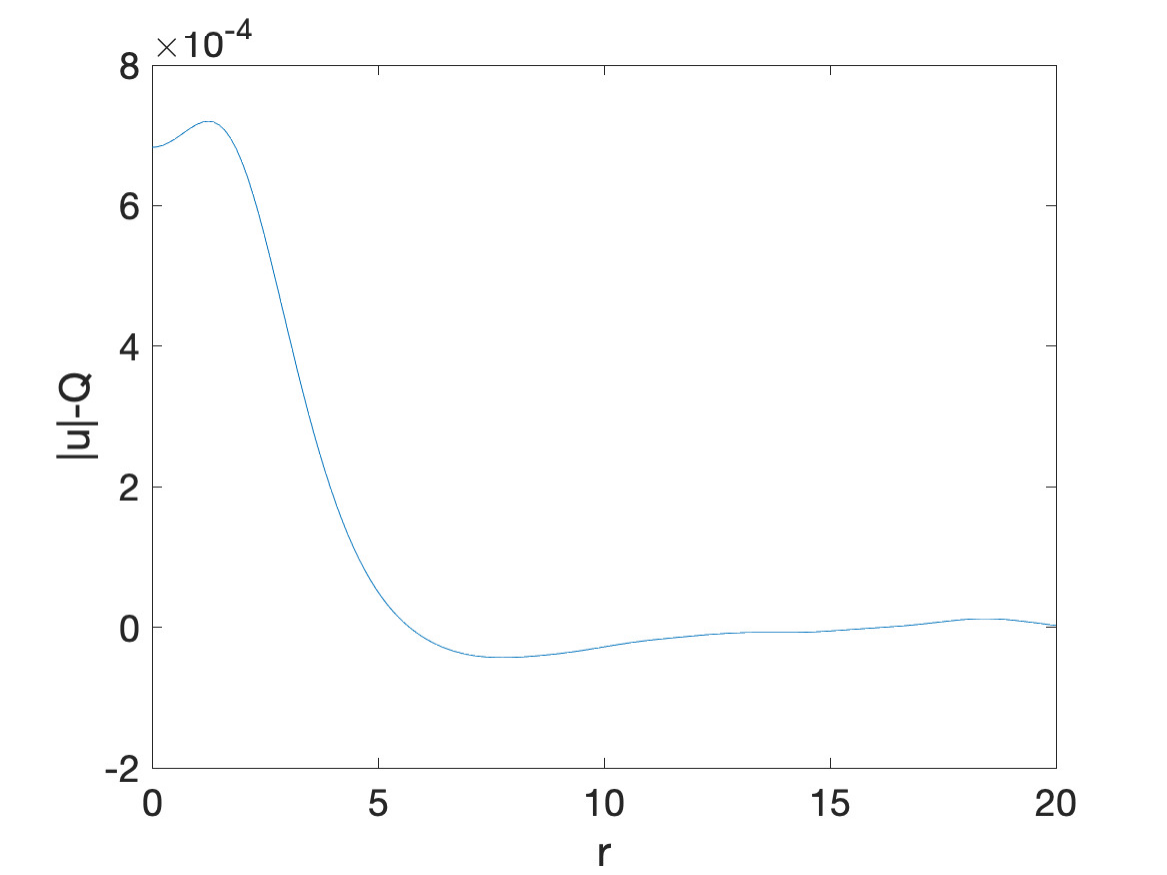}
 \caption{The 2D solution $u$ to \eqref{eq:nls} for initial data 
 (\ref{ini}) with $\omega=0.1$ and $\lambda =1.001$. On the left the $L^{\infty}$-norm as a function of time. On the right the difference between $|u|$ and $Q_{\rm \omega=0.1}$ at the final time $t_{\rm f}=20$.}
 \label{NL35_d2solom011001}
\end{figure}

As a second case, we consider the same 2D initial data \eqref{ini}, but with $\lambda=1.001$. In Fig.~\ref{NL35_d2solom011001} we again 
show the $L^{\infty}$-norm of the solution as a function of time. 
Similarly as before, a final state is reached and its maximum is again
found to be very close to the unperturbed ground state $Q_{\omega=0.1}$. In both cases, we find that the difference is largest for $r$ close to the origin.

\begin{figure}[htb!]
  \includegraphics[width=0.49\textwidth]{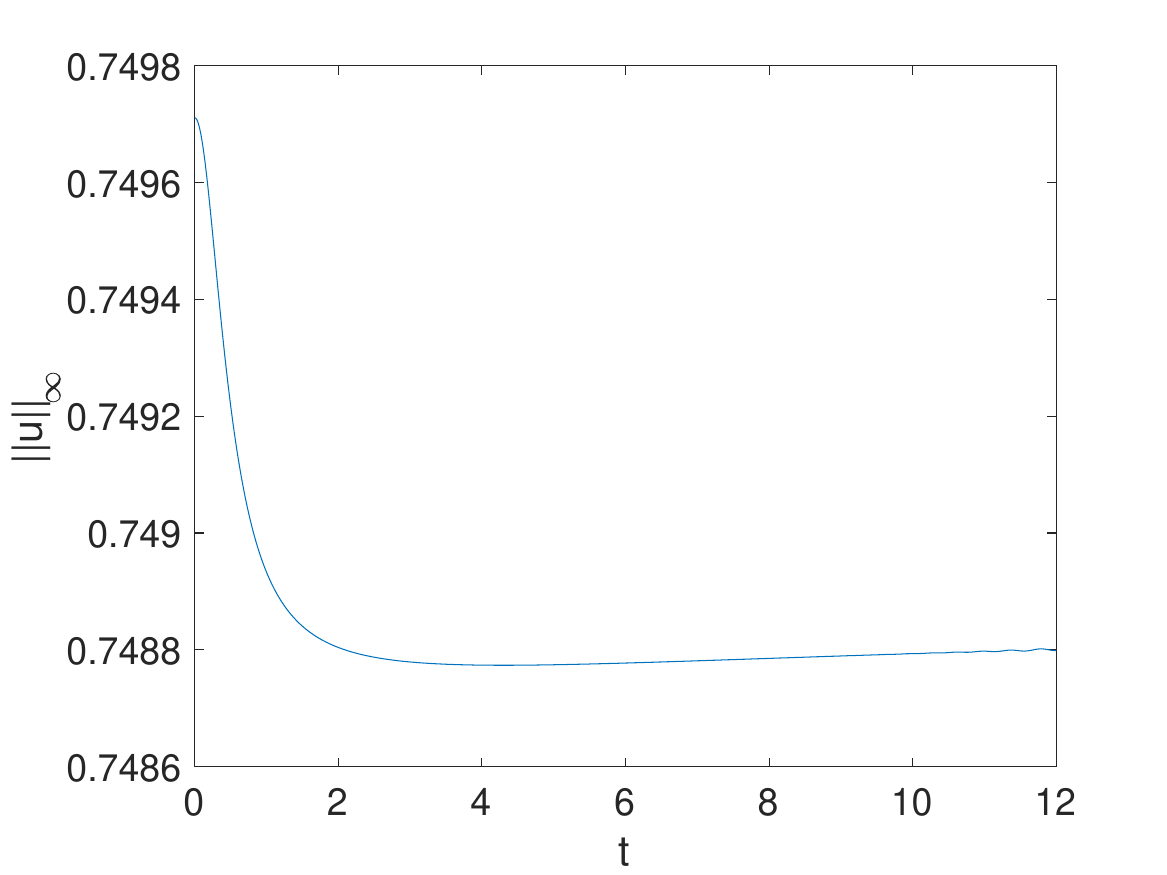}
  \includegraphics[width=0.49\textwidth]{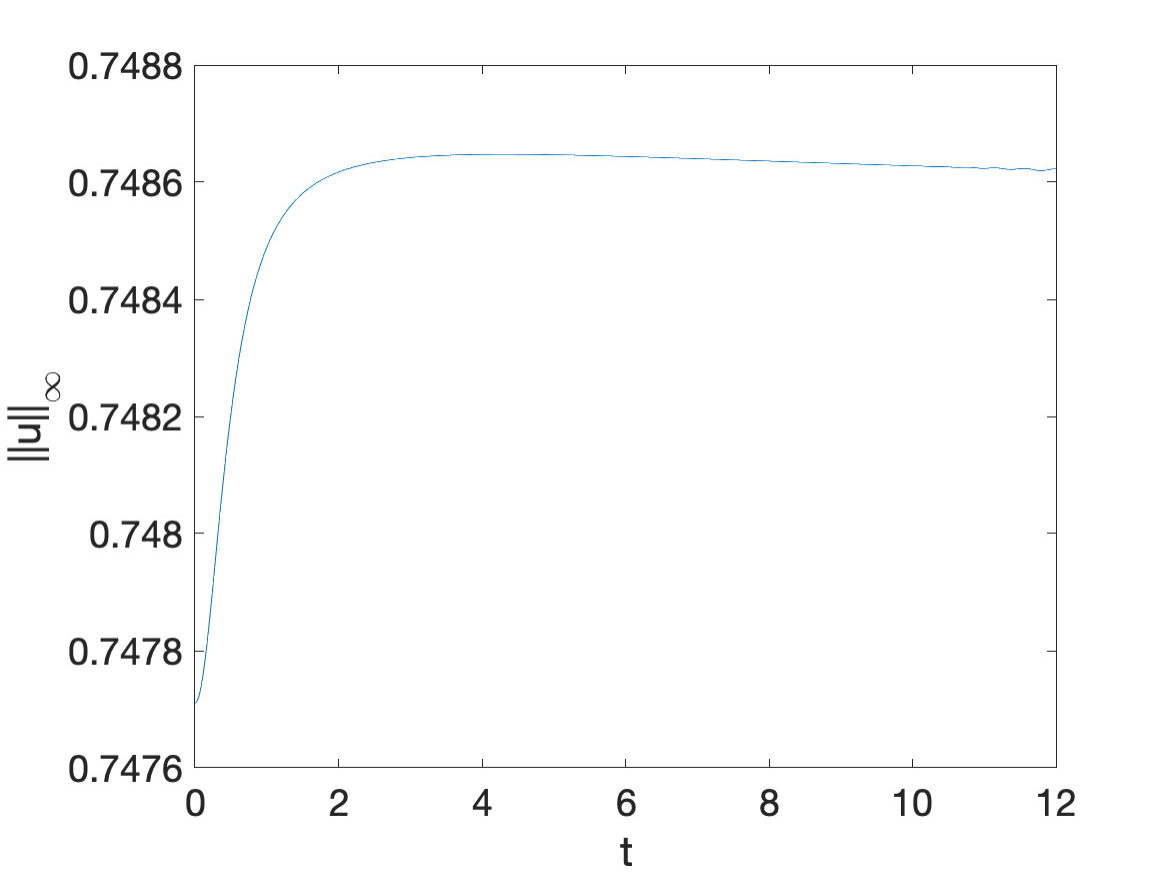}
 \caption{$L^{\infty}$-norm as a function of time for the 2D solution $u$ to \eqref{eq:nls} with initial data 
 \eqref{inibis} and $\omega=0.1$. On the left the case with ``$+$" sign. On the right the case with ``$-$" sign.}
 \label{NL35_d2solom01gauss}
\end{figure}

In order to illustrate that the qualitative picture found before is not due to our specific choice of 
perturbations, we shall also consider ground states perturbed by a small Gaussian-like perturbation, i.e.
\begin{equation}
	u_0(x) =  Q_{\omega=0.1}(x)\pm \lambda e^{-|x|^{2}}, \quad \lambda = 0.001.
	\label{inibis}
\end{equation}
Note that we only consider smooth perturbations in 
this paper in order to allow for spectral accuracy in the radial 
coordinate, i.e., an exponential decrease of the numerical error with 
the number of collocation points.
In Fig.~\ref{NL35_d2solom01gauss} we show the behavior in time of the respective $L^{\infty}$-norms for the two choices $\pm \lambda$. In 
both situations the difference between $|u|$ at $t_{\rm f}=20$ and $Q_{\rm \omega =0.1}$ is found to be of the order $\mathcal O(10^{-4})$. Moreover, 
the error (not depicted here for the sake of readability) is again found to be largest close to the origin.

\begin{figure}[htb!]
  \includegraphics[width=0.49\textwidth]{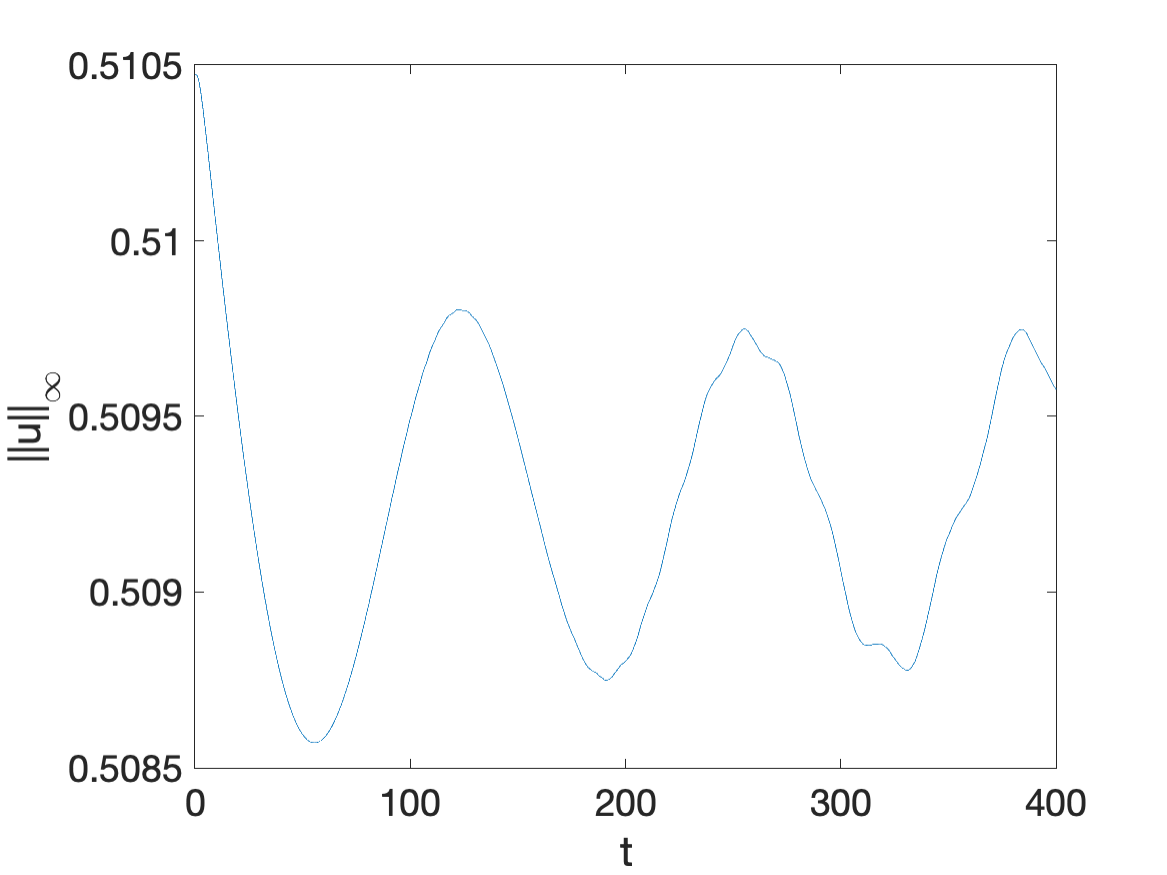}
  \includegraphics[width=0.49\textwidth]{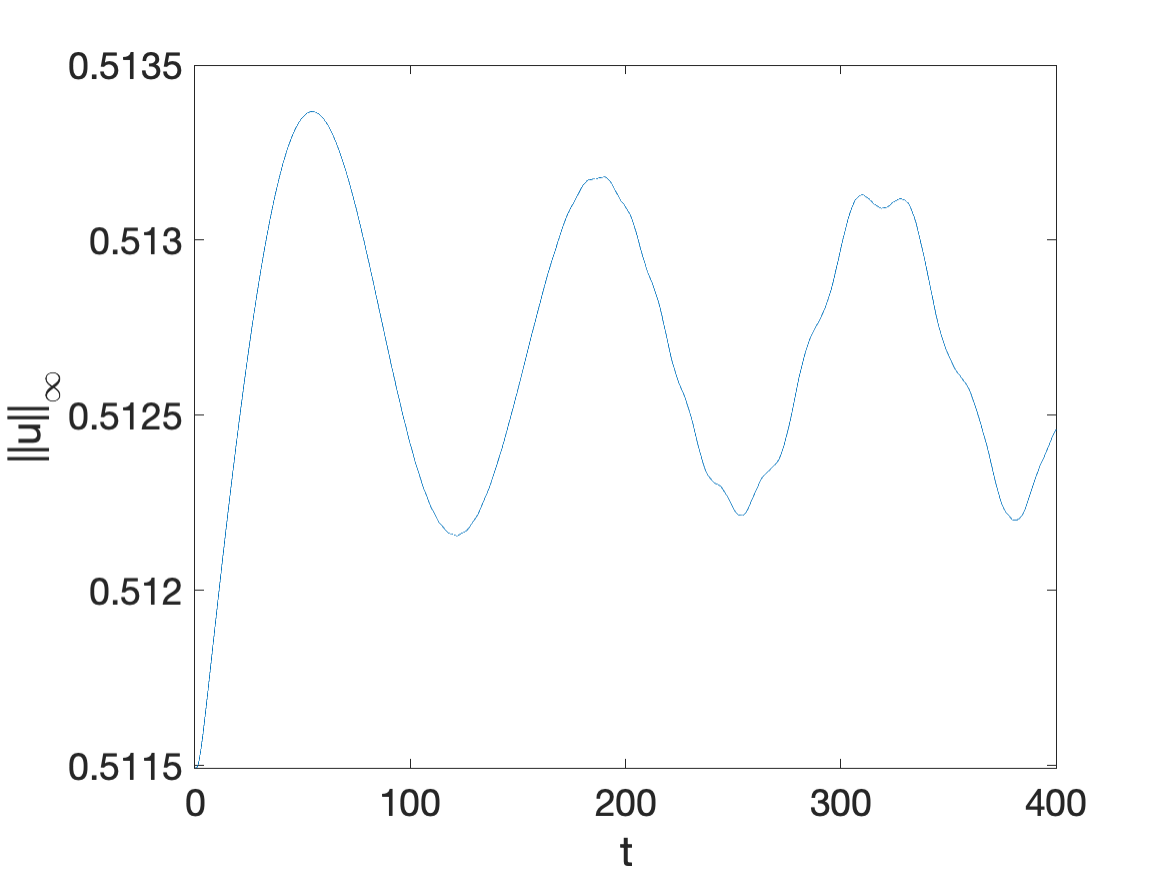}
 \caption{$L^{\infty}$-norm as a function of time for the 2D solution $u$ to \eqref{eq:nls} with initial data 
 \eqref{inibis} and $\omega=0.05$. On the left the case 
  with the ``$+$" sign. On the right the case with the ``$-$" sign.}
 \label{NL35_d2solom005}
\end{figure}

If similar perturbations are applied to other ground states $Q_\omega$, 
the resulting solution $u$ behaves qualitatively similarly. 
Our numerical tests therefore support
Conjecture~\ref{conj:2d}. However, we also find that the smaller the choice of $\omega\in (0, \tfrac{3}{16})$, the longer 
it takes for the solution $u$ to reach its final state. In fact, for small enough $\omega$, damping effects within 
the time-oscillations of $|u|$ become almost invisible, even if one computes up to much larger times 
$t_{\rm f}=400$, see Fig.~\ref{NL35_d2solom005}. 


\section{(In-)stability of action  ground states in 3D}\label{sec:stab3D}

\subsection{Stable branch} In this section, we shall study the question of (in-)stability of cubic-quintic ground states 
in dimension $d=3$. In view of Fig. \ref{NL35_d3solmass}, we expect ground states $Q_\omega$ 
with $\omega > \omega_{\rm crit}\approx 0.026$ to be orbitally stable. That this is indeed the case, is strongly suggested 
by our numerical results below.

\begin{figure}[htb!]
  \includegraphics[width=0.49\textwidth]{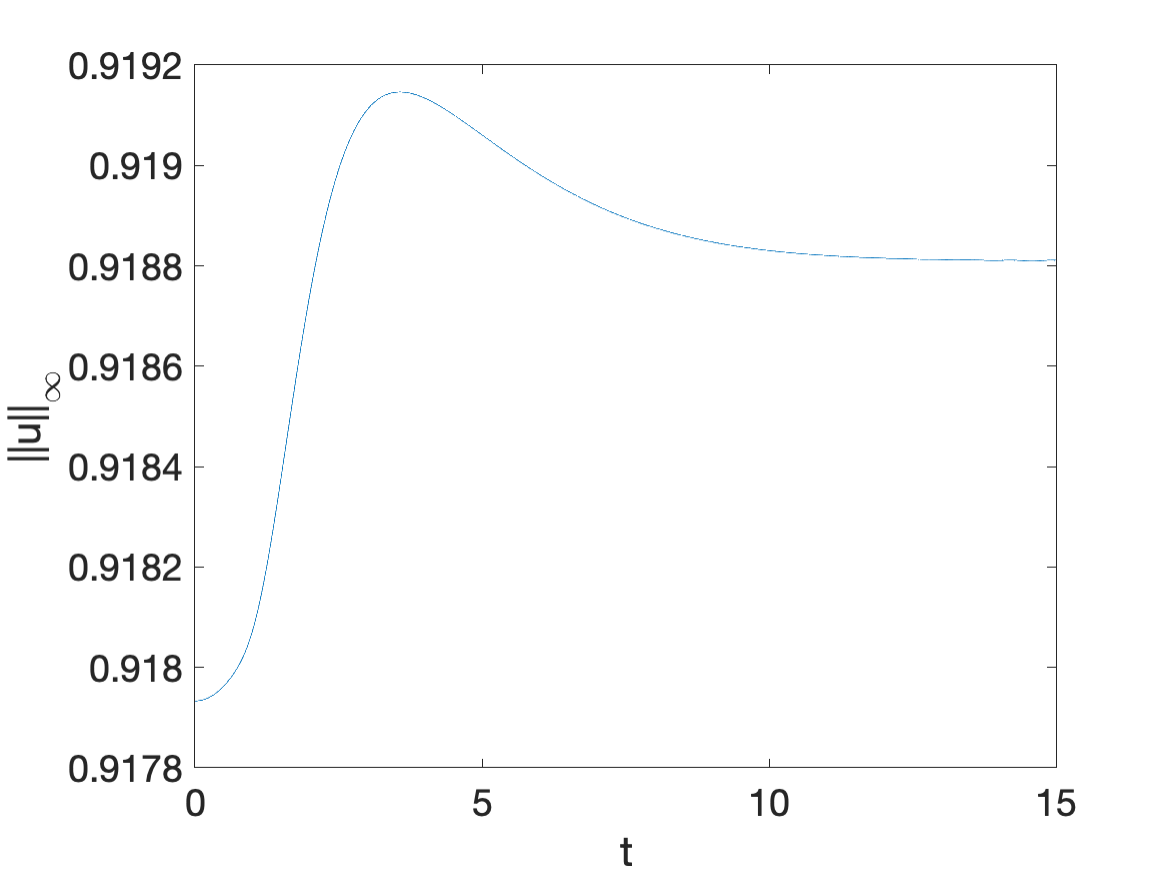}
  \includegraphics[width=0.49\textwidth]{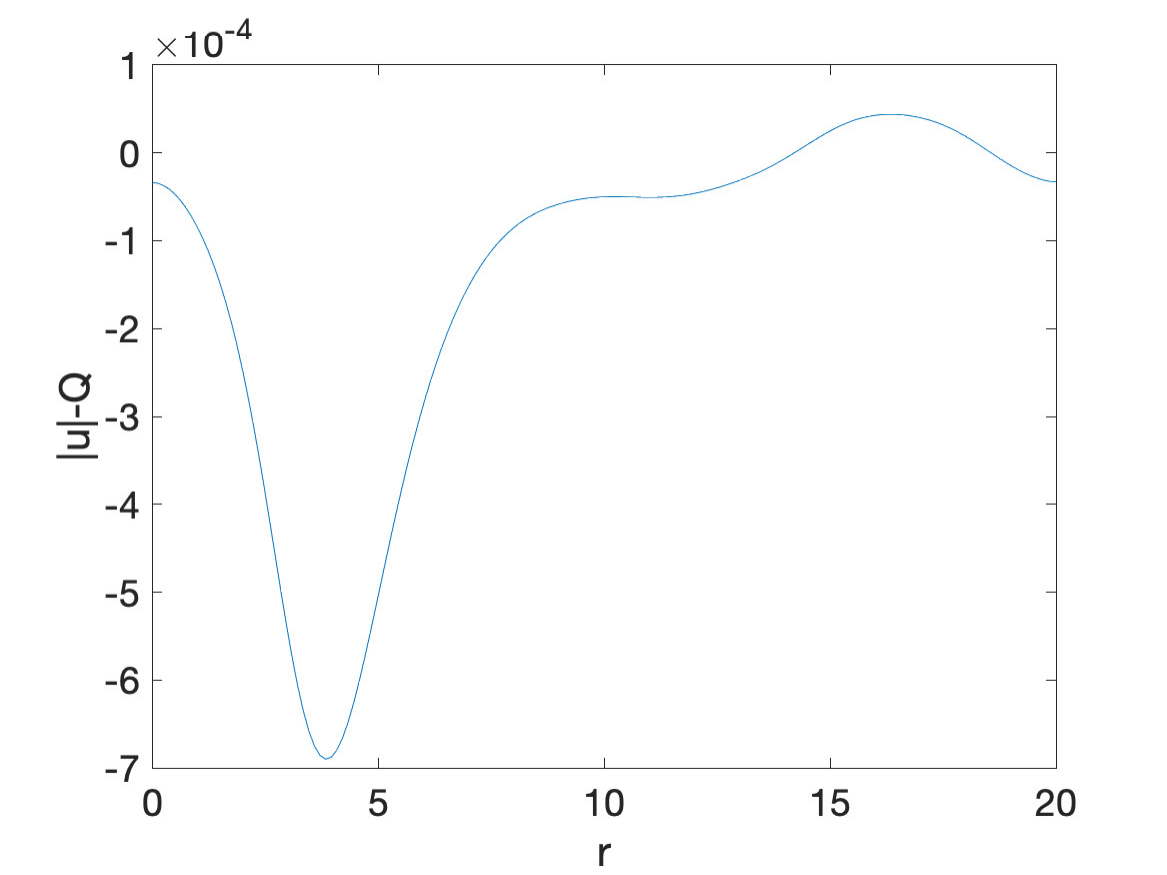}
 \caption{The 3D solution $u$ to \eqref{eq:nls} for initial data 
 (\ref{ini}) with $\omega=0.1$ and $\lambda =0.99$. On the left the $L^{\infty}$-norm as a function of time. On the right 
 the difference between $|u|$ and $Q_{\rm \omega=0.1}$ at the final time $t_{\rm f}=15$.}
 \label{NL35_d3solom01}
\end{figure}

\begin{figure}[htb!]
  \includegraphics[width=0.49\textwidth]{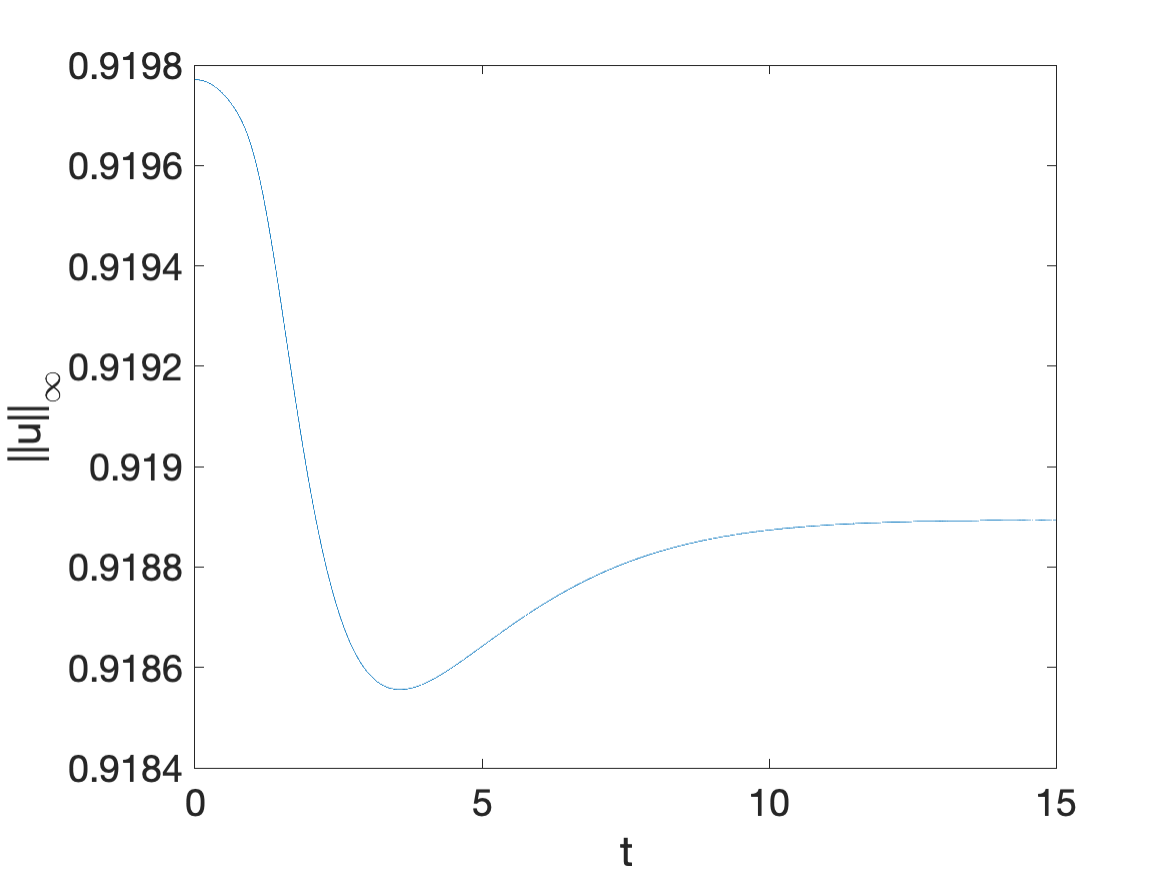}
  \includegraphics[width=0.49\textwidth]{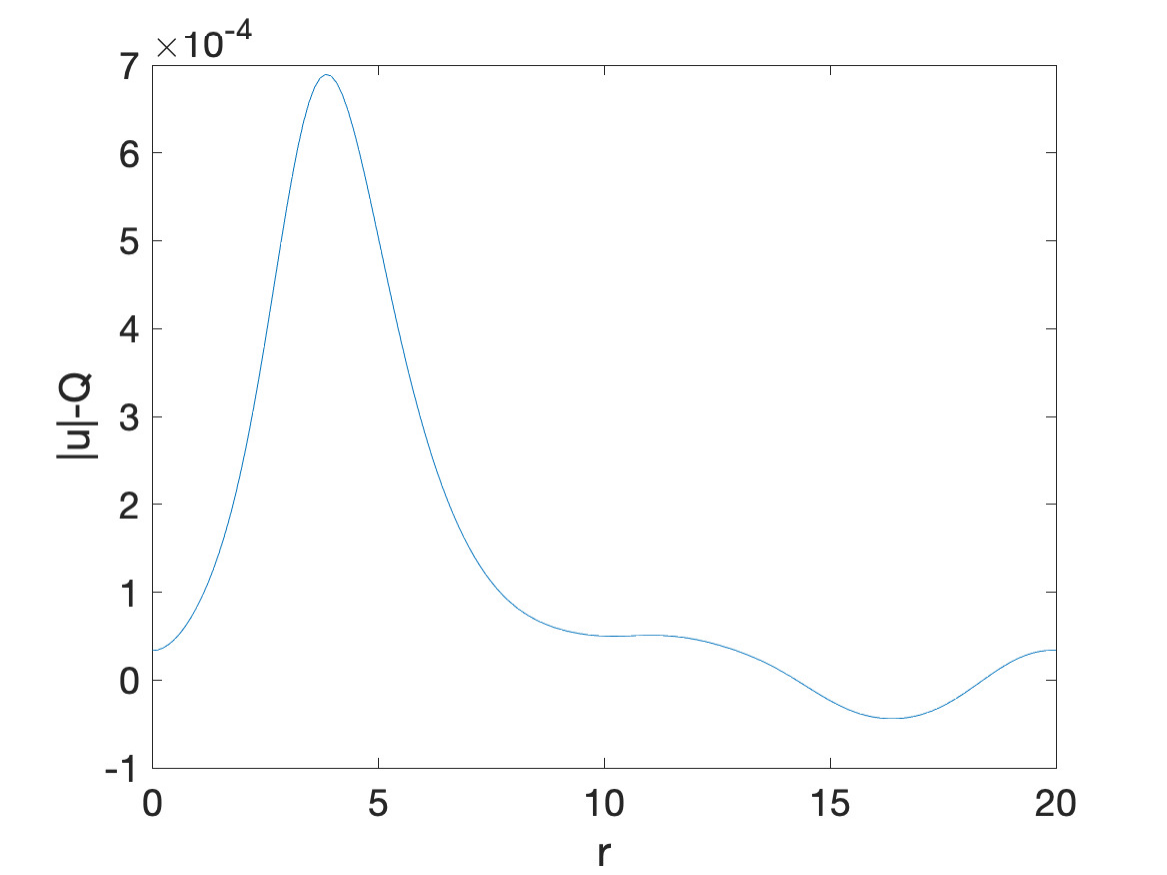}
 \caption{The 3D solution $u$ to \eqref{eq:nls} for initial data 
 (\ref{ini}) with $\omega=0.1$ and $\lambda =1.001$. On the left the $L^{\infty}$-norm as a function of time. On the right 
 the difference between $|u|$ and $Q_{\rm \omega=0.1}$ at the final time $t_{\rm f}=15$.}
 \label{NL35_d3solom011001}
\end{figure}

To this end, we first consider multiplicative perturbations of $Q_{\rm \omega}$ on the stable branch: 
In Figures \ref{NL35_d3solom01} and  \ref{NL35_d3solom011001} we study the time-evolution of \eqref{eq:nls} with 
initial data of the form \eqref{ini}. 
On the left of Fig.~\ref{NL35_d3solom01} we show the $L^{\infty}$-norm of the solution $u$ obtained in the case $\omega =0.1$ and 
$\lambda=0.99$. On the right of the same figure, we show the difference between the unperturbed ground state $Q_{\rm \omega = 0.1}$ and $|u|$ at the 
final time $t_{\rm f} = 15$. It can be 
seen that the $L^{\infty}$ norm settles on a nearly constant value as $t\to 15$.

In Fig.~\ref{NL35_d3solom011001} we study the analogous situation with $\lambda=1.001$: Again, the (absolute value of the) solution $u$ 
seems to settle around $t_{\rm f} = 15$ on the stable unperturbed ground state $Q_{\rm 0.1}$. In both cases, the error between $|u|$ and 
$Q_{\rm \omega=0.1}$ is again found to be of the order $\mathcal O(10^{-4})$.


\subsection{Unstable branch} 

\begin{figure}[htb!]
  \includegraphics[width=0.49\textwidth]{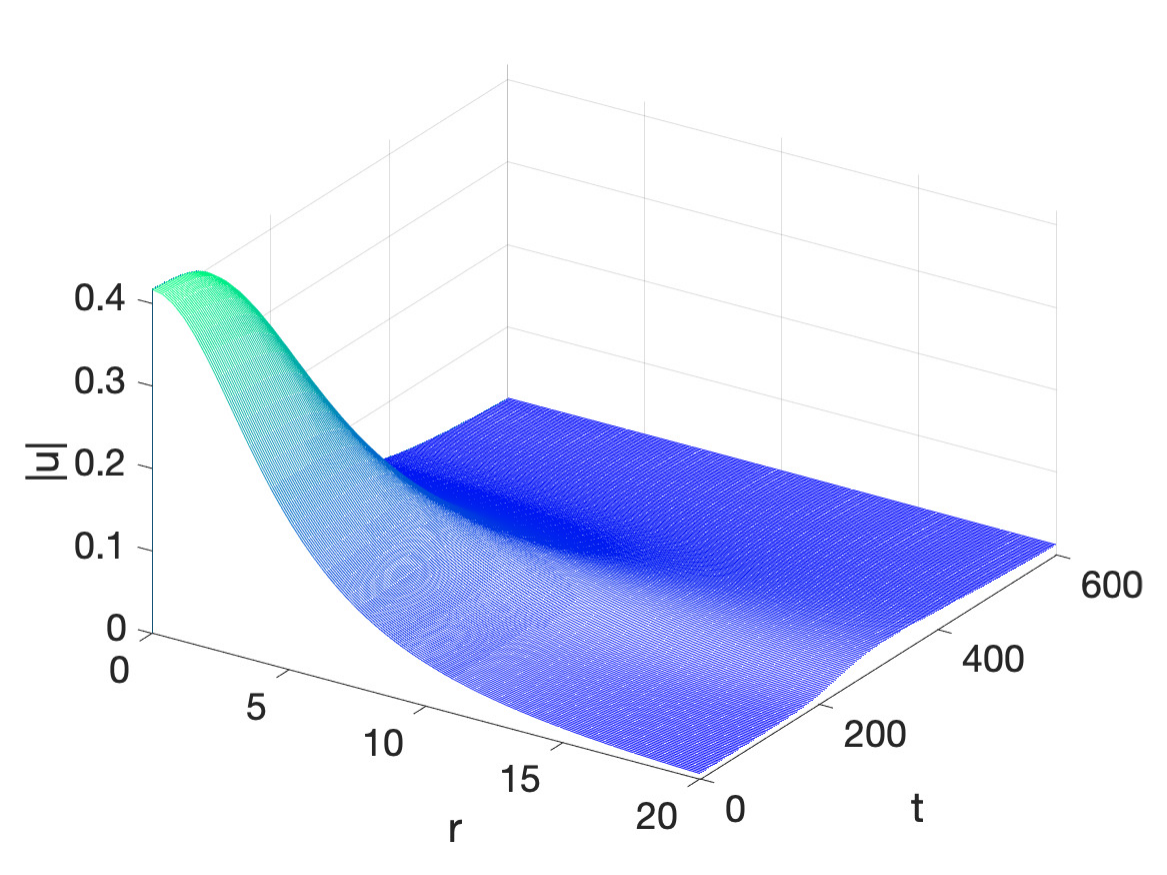}
  \includegraphics[width=0.49\textwidth]{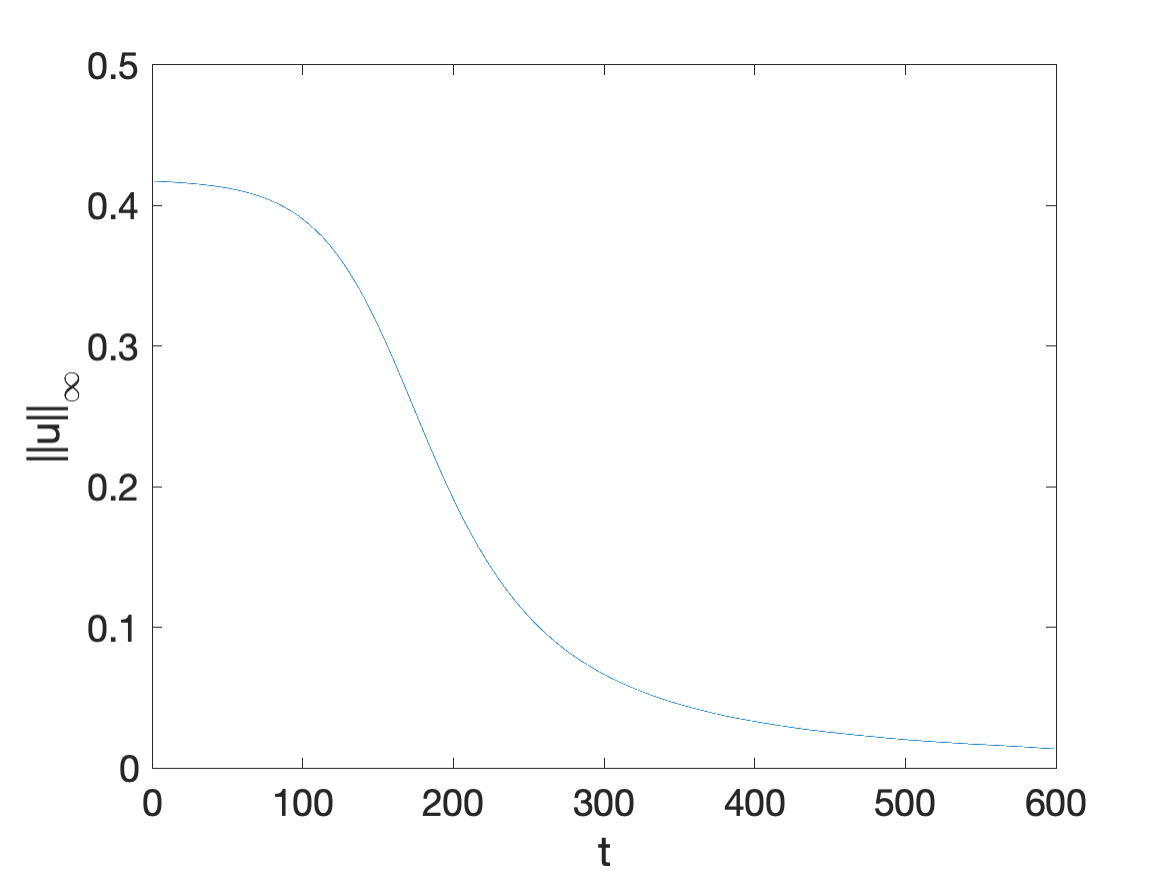}
 \caption{Left: The absolute value of the solution $u$ to \eqref{eq:nls} in dimension $d=3$ 
 obtained from initial data \eqref{ini} with $\omega=0.01$ and $\lambda=0.999$. Right: The $L^{\infty}$-norm of $u$ as a function of time.}
 \label{NL35_d3solom002}
\end{figure}

The situation dramatically changes if we consider perturbations of ground state solutions on the unstable branch, i.e. perturbations of $Q_{\omega}$ with $\omega<\omega_{\rm crit}\approx 0.026$:

In Fig.~\ref{NL35_d3solom002} we show  
the solution $u$ to \eqref{eq:nls} obtained from initial data \eqref{ini} with $\omega=0.01$ and $\lambda=0.999$. 
Note that this implies $M(u_0)< M(Q_{\omega})$. 
The solution is seen to be purely dispersive which is also confirmed by 
the $L^{\infty}$-norm of the solution as a function of time (depicted on the 
right of the same figure). In fact, we did not 
discover any stable structure within the time-evolution even if we let the numerical code run 
for longer times. 

\begin{figure}[htb!]
  \includegraphics[width=0.7\textwidth]{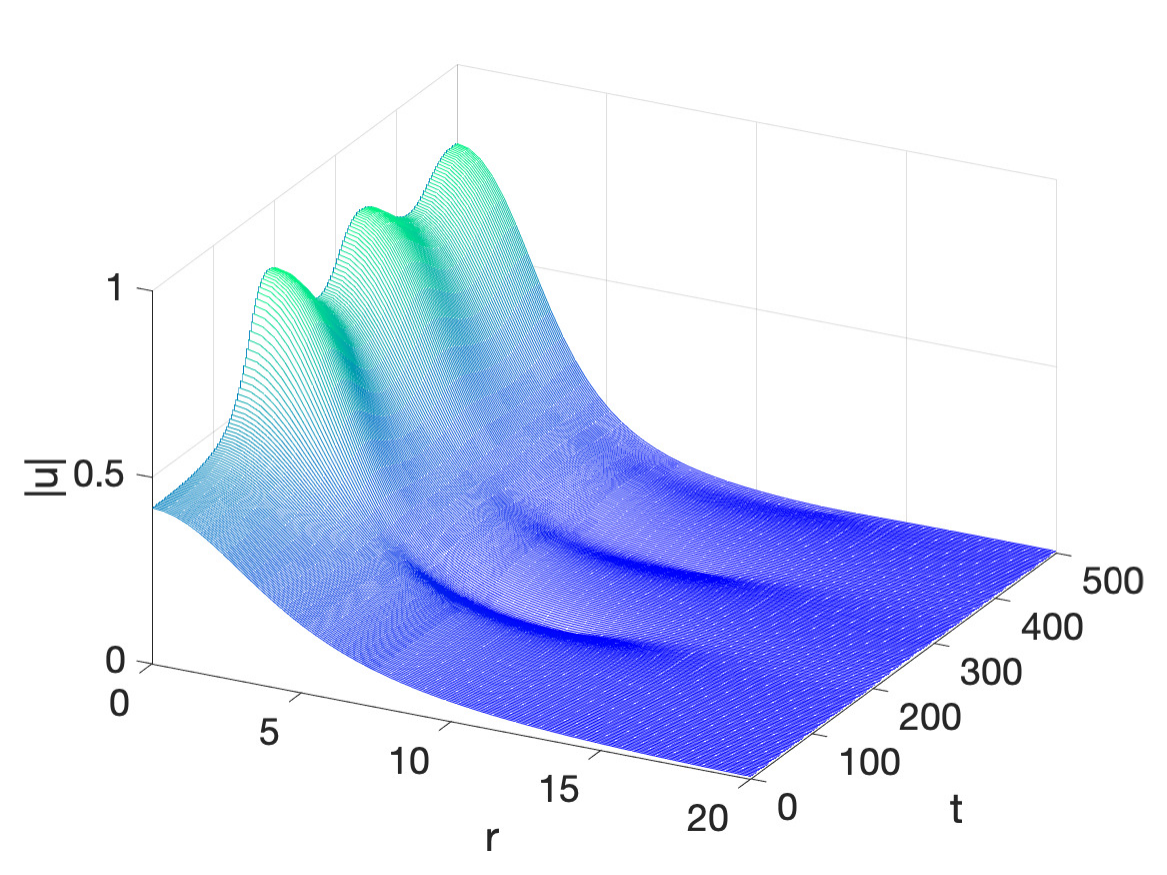}
 \caption{Solution to \eqref{eq:nls} in dimension $d=3$ 
 for initial data (\ref{ini}) with $\omega=0.01$ and $\lambda=1.001$.}
 \label{NL35_d3solom002101}
\end{figure}

If we consider the same ground state as before, but instead choose $\lambda=1.001$, we find a different kind of instability. 
Now the solution $u$ shows oscillations of high amplitude, see Figure \ref{NL35_d3solom002101}. 
\begin{figure}[htb!]
  \includegraphics[width=0.49\textwidth]{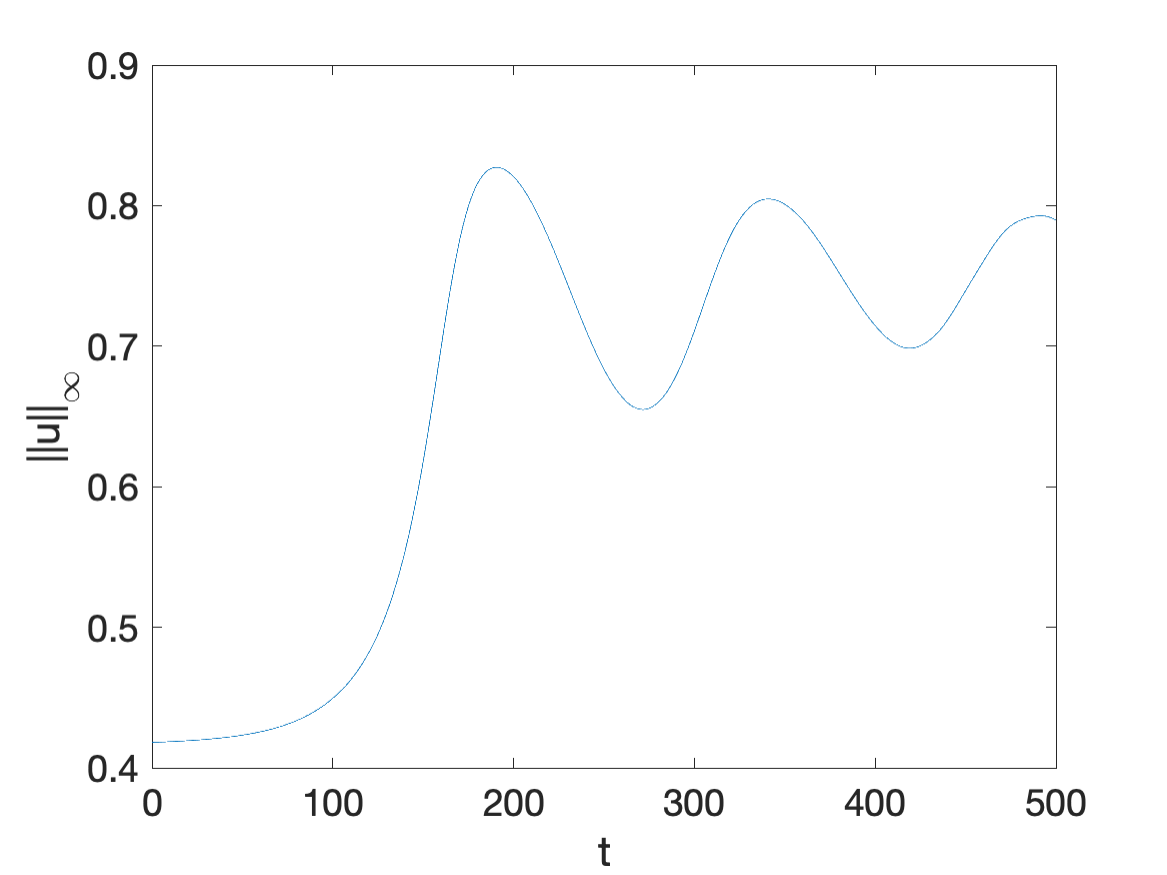}
  \includegraphics[width=0.49\textwidth]{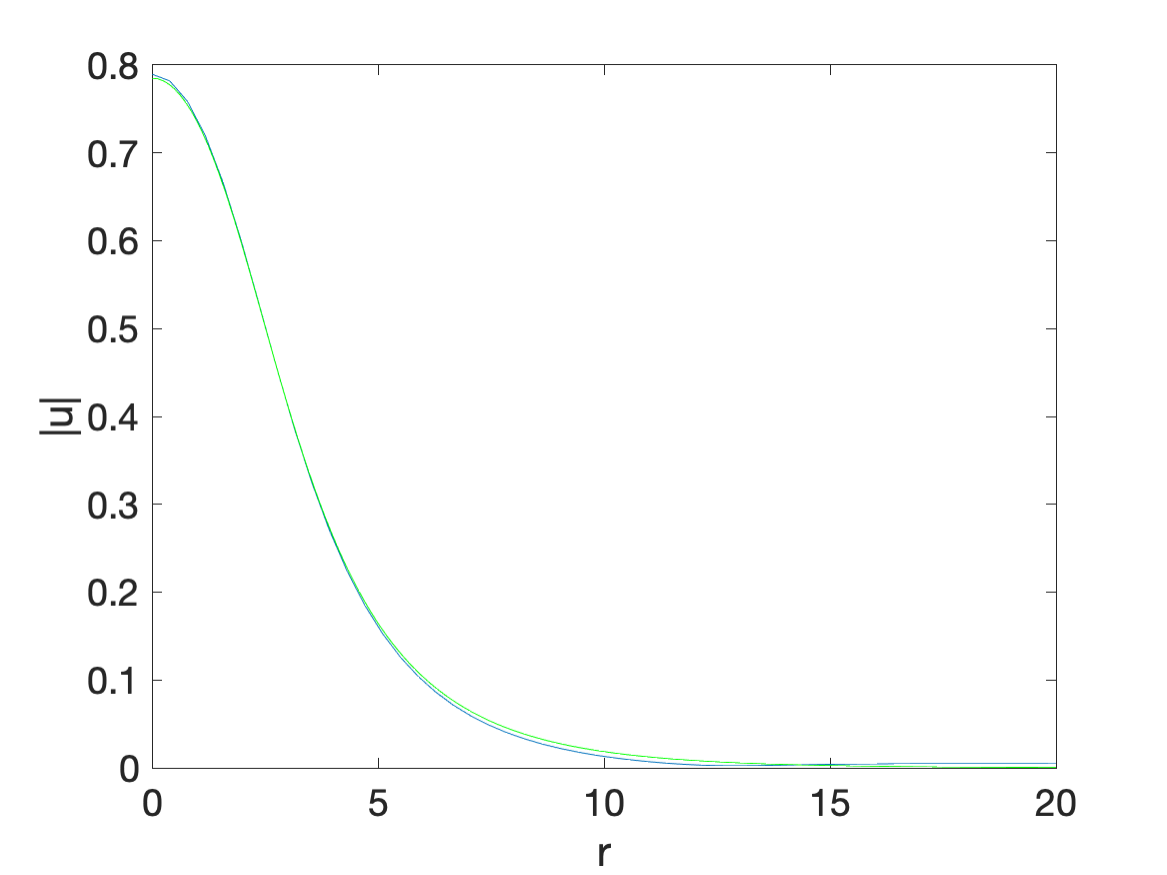}
 \caption{Solution to \eqref{eq:nls} in dimension $d=3$
 for initial data (\ref{ini}) with $\omega=0.01$ and $\lambda=1.001$: On the 
 left the $L^{\infty}$-norm of $u$ as a function of time. On the right $|u|$ at the final 
 time (in blue) together with the ground state $Q_{\omega=0.047}$ (in green).}
 \label{NL35_d3solom002101inf}
\end{figure}

These oscillations are even more visible in the $L^{\infty}$-norm of the solution, as depicted 
on the left of Fig.~\ref{NL35_d3solom002101inf}. 
One can see that early on the norm is growing 
strongly but then it appears to show damped 
oscillations around some final state. We conjecture that the latter 
corresponds to another ground state on the stable branch.  To this end, we compare the maximum of $|u|$, obtained at the final time $t_{\rm f}=500$, with 
the $L^\infty$-norms in our library of previously computed action  ground states $Q_{\rm \omega}$, cf. Fig.~\ref{NL35_d3sol}. 
Indeed we find good agreement of $|u|$, when compared to $Q_{\omega}$ with $\omega=0.047>\omega_{\rm crit}$, 
see the right of Fig.~\ref{NL35_d3solom002101inf}. 
Thus perturbations of unstable ground states where $M(u_0)>M(Q_\omega)$, seem to result in solutions which eventually settle on 
another, stable ground state as $t\to+ \infty$. 
Note, however, that $M(Q_{\omega=0.01})\approx 79.44$ 
while $M(Q_{\omega=0.047})\approx 77.05$. If the final state had the same 
mass as the unperturbed initial state, this would correspond to an 
$\omega\approx 0.0495$. This shows that a certain amount of  mass is lost through radiation.


\subsection{Other kinds of perturbations}\label{sec:other-perturb}
The results described above are not due to our specific choice of perturbations. To show this, we 
consider initial data
\begin{equation}\label{ini2}
	u_{0, \pm}(x) = Q_{\omega}(x)\pm \lambda e^{-\left(|x|-|x_0|\right)^{2}}, \quad \lambda =0.001.
\end{equation}

\begin{figure}[htb!]
  \includegraphics[width=0.49\textwidth]{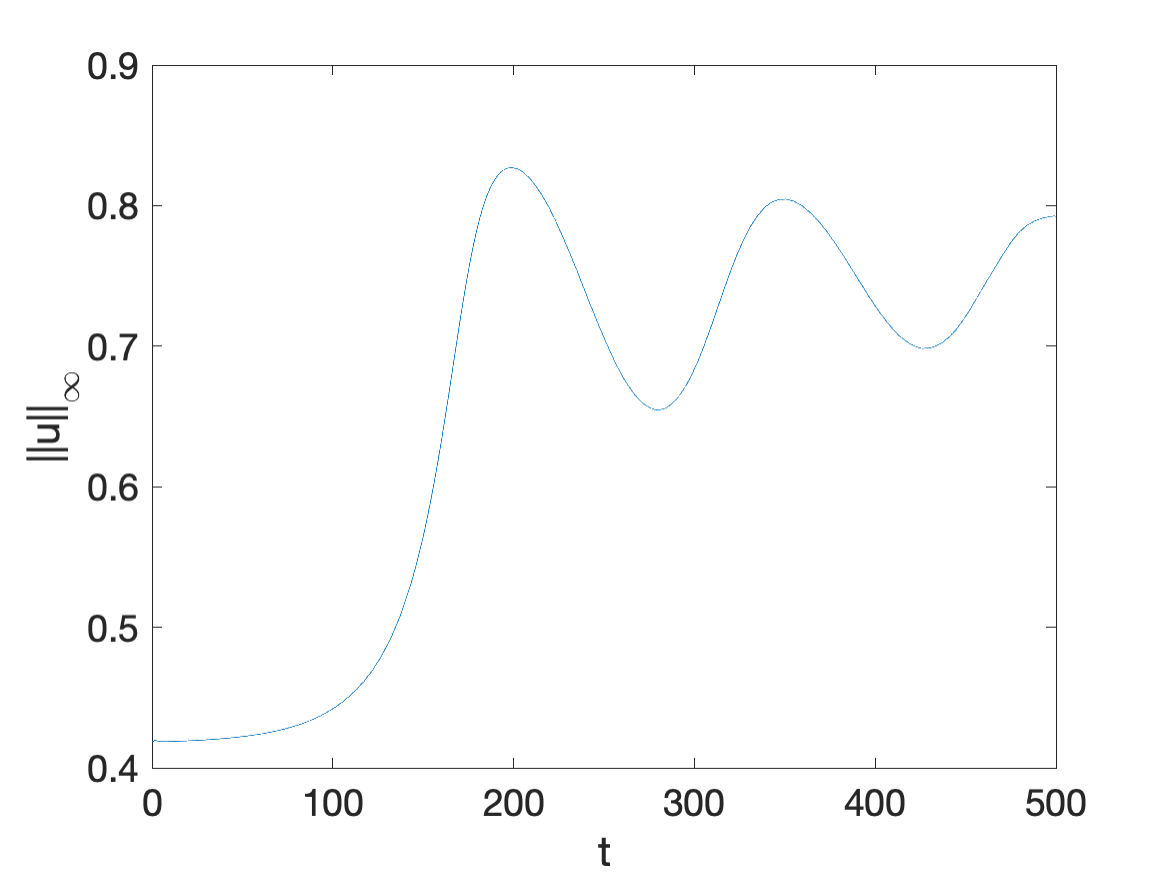}
  \includegraphics[width=0.49\textwidth]{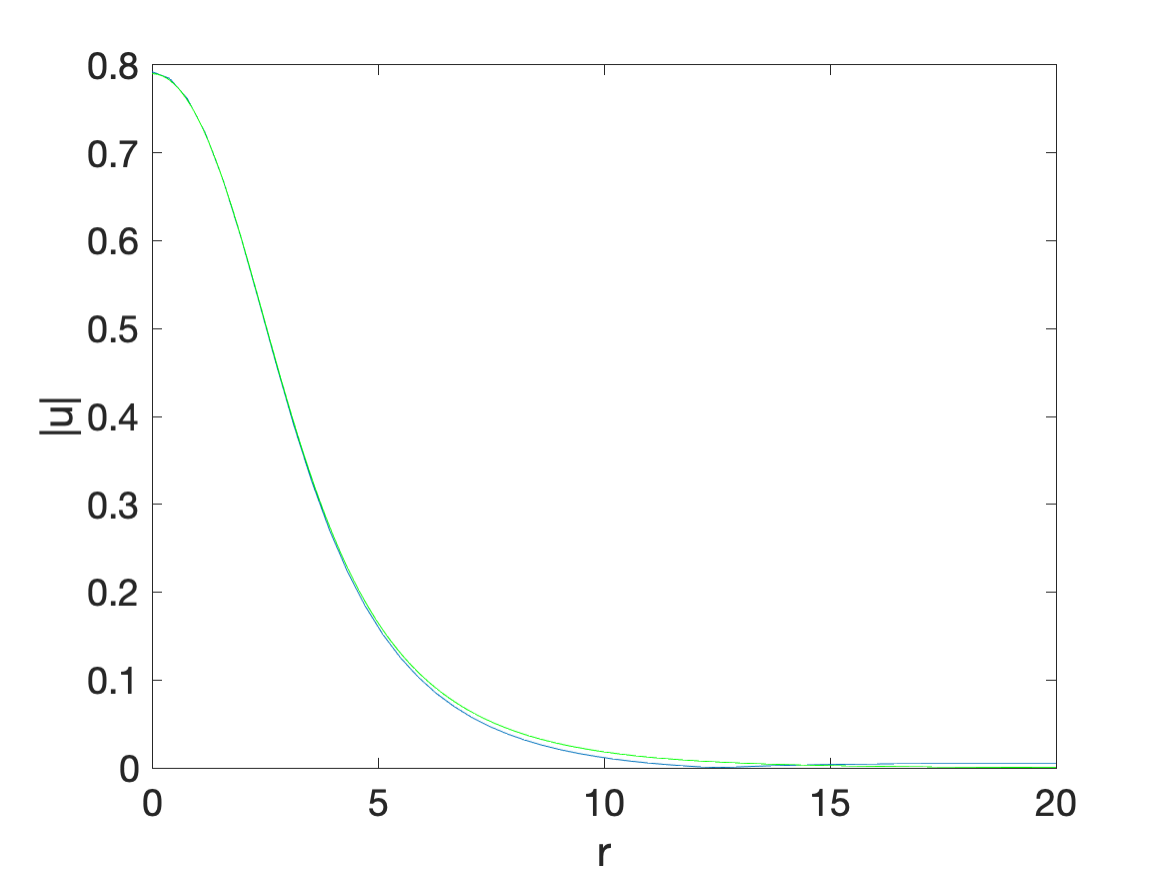}
 \caption{Solution to \eqref{eq:nls} in dimension $d=3$
 for initial data of the form \eqref{ini2} with $|x_0|=1$, and $\omega=0.01$: On the left the $L^{\infty}$-norm of $u$ as a function of time. 
 On the right $|u|$ at the final 
 time (in blue) together with the ground state $Q_{\omega=0.048}$ (in green).}
 \label{NLS35_d3solom001p001sgaussinf}
\end{figure}

For both $|x_0|=0$ and $|x_0| \not =0$, the solution in the case with the ``$+$" sign looks very similar to the one depicted in 
Fig.~\ref{NL35_d3solom002101inf}. This fact becomes particularly clear when one compares the 
time-evolution of the $L^{\infty}$-norm of $u$ depicted in Fig.~\ref{NLS35_d3solom001p001sgaussinf}, with the one from Fig.~\ref{NL35_d3solom002101inf}.

By comparing the maximum of $|u|$ found at the
final time $t_{\rm f}=600$ with  
the $L^\infty$-norm of a stable ground state, we find good agreement
with $Q_{\omega=0.048}$. The latter has mass $M(Q_{\omega=0.048})\approx 77.95$. Unfortunately, we are unable 
to clearly decide whether the final state is closer to $Q_{\omega=0.048}$ than to $Q_{\omega=0.047}$.

\begin{figure}[htb!]
  \includegraphics[width=0.49\textwidth]{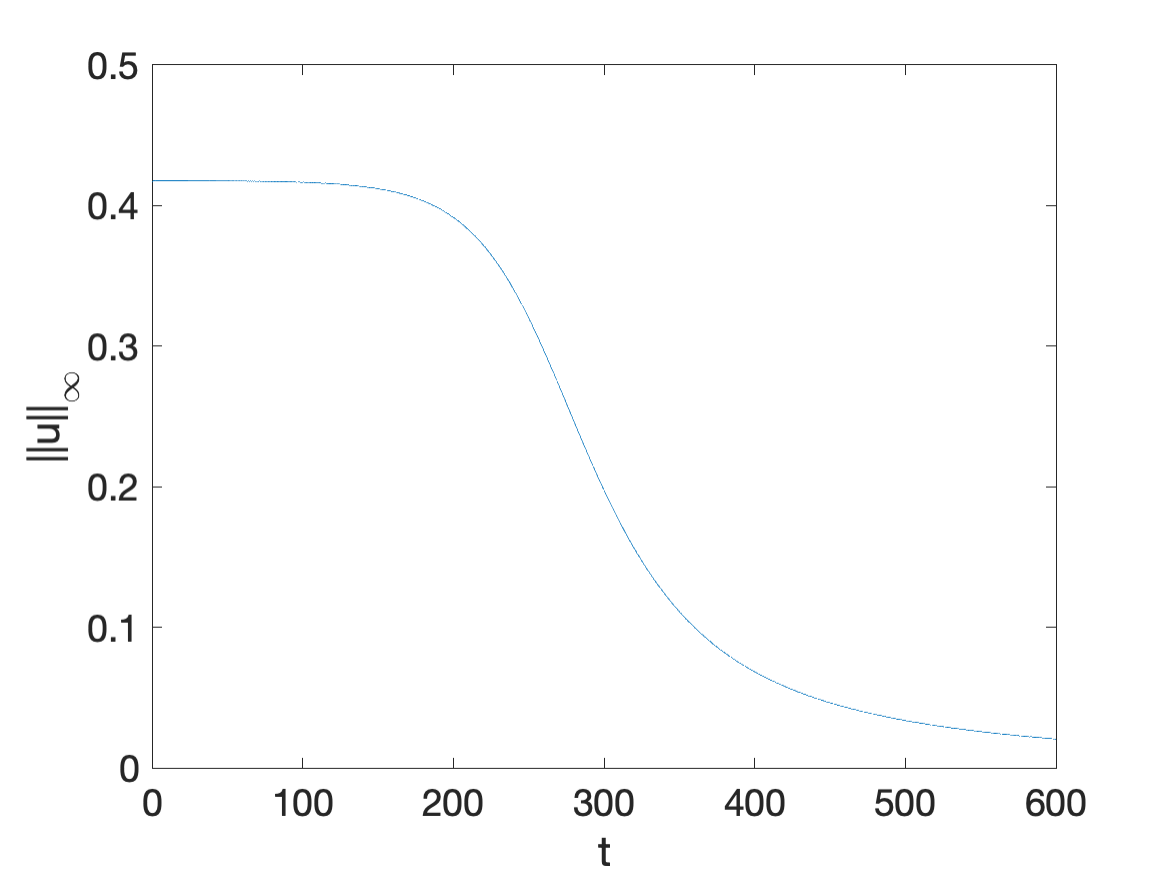}
  \includegraphics[width=0.49\textwidth]{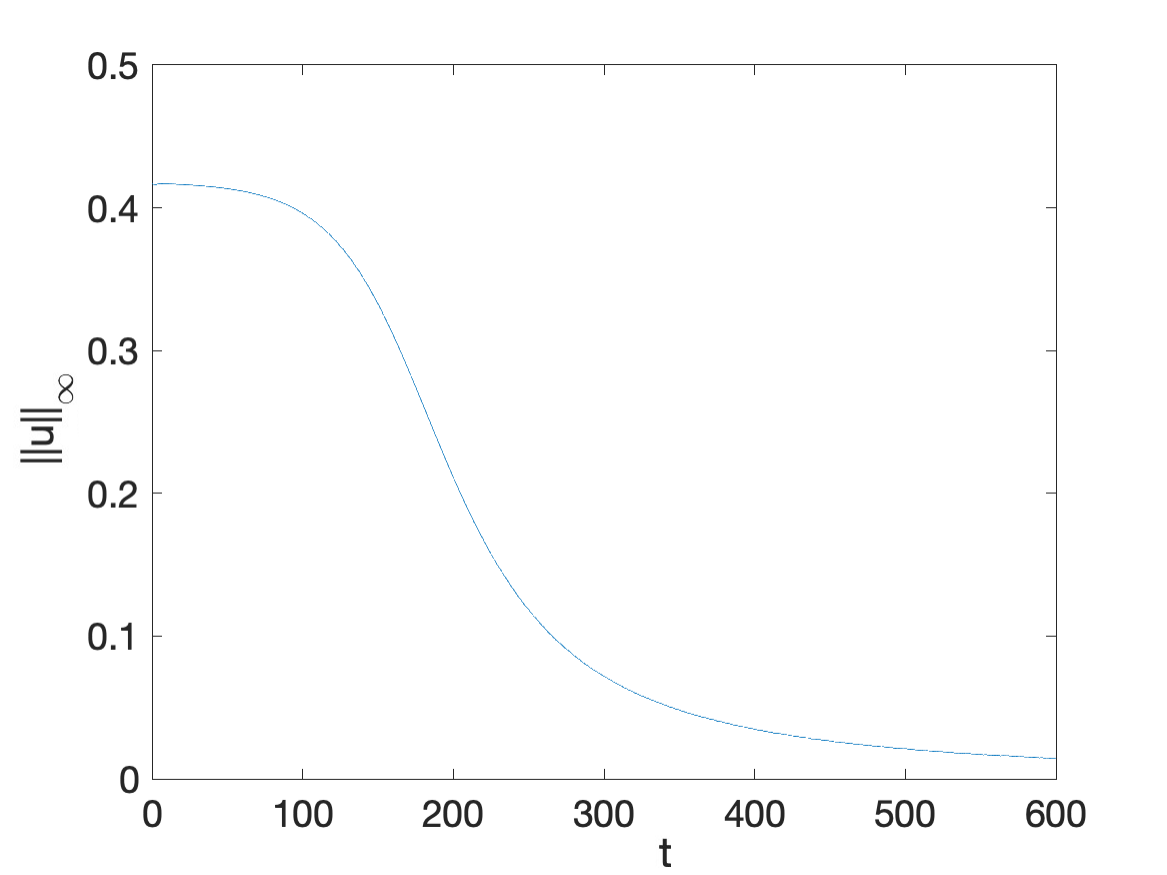}
  \caption{The $L^{\infty}$-norm of the solution $u$ in dimension $d=3$ obtained from initial data  \eqref{ini2} with the ``$-$" sign. 
  On the left the case with $x_0=0$ and on the right the one with $x_0=1$.}
  \label{NLS35_d3solom001m001gaussinf}
\end{figure}

In the case of initial data \eqref{ini2} with the ``$-$" sign, we 
again find that the solution is completely dispersed, see Fig.~\ref{NLS35_d3solom001m001gaussinf}. 
This is consistent with our earlier findings above which indicate that 
perturbation with $M(u_0)<M(Q_\omega)$ lead to purely dispersive solutions.

\begin{figure}[htb!]
  \includegraphics[width=0.49\textwidth]{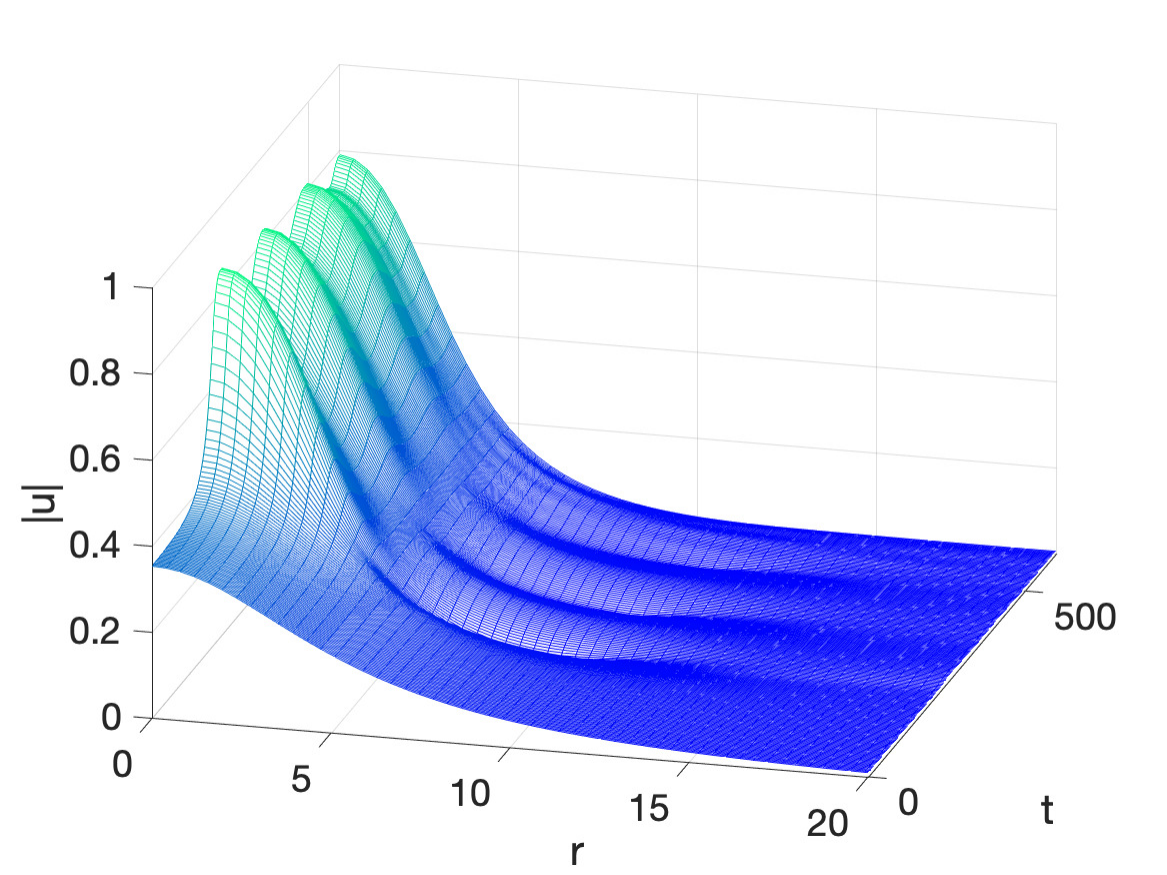}
  \includegraphics[width=0.49\textwidth]{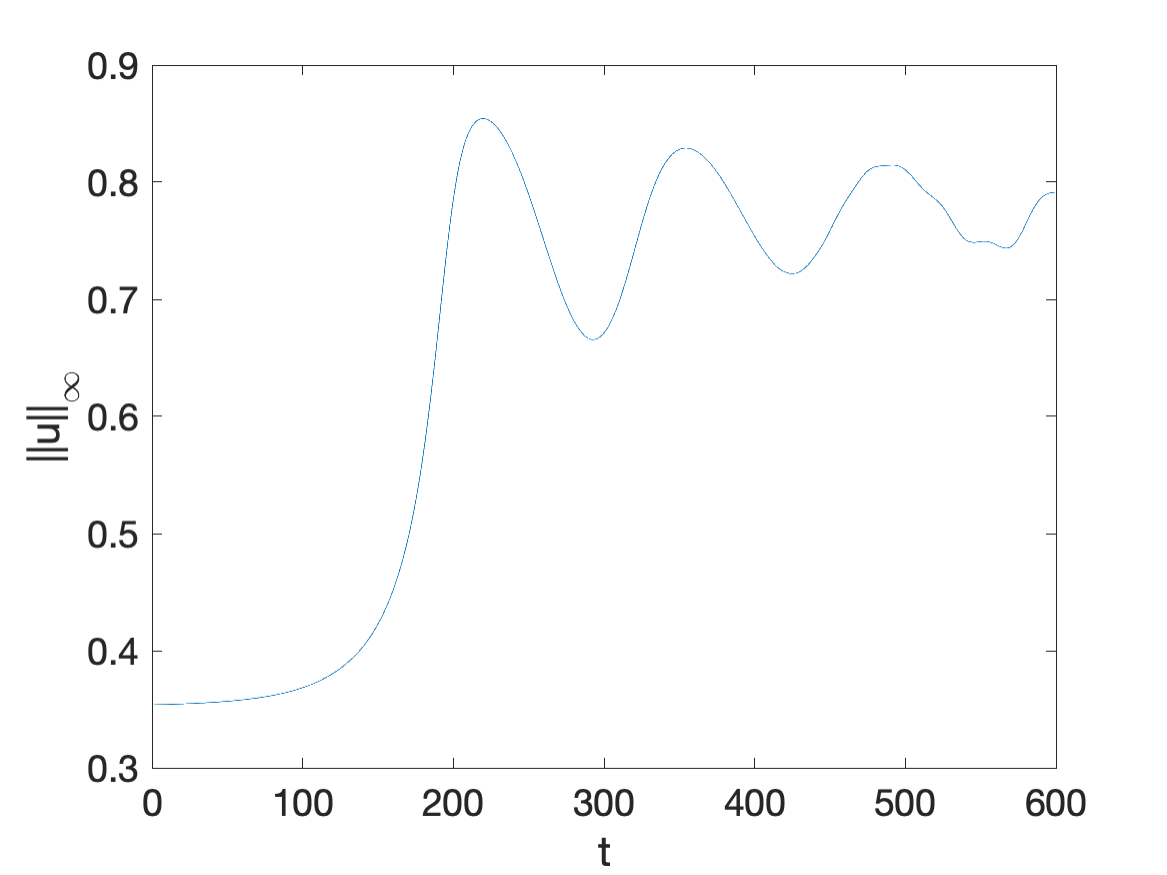}
  \caption{Left: Solution to \eqref{eq:nls} in dimension $d=3$ for initial data 
  (\ref{ini2}) with $x_0=0$ and $\omega=0.007$. Right: 
  The $L^{\infty}$-norm of $u$ as a function of time.}
  \label{NLS35_d3solom001m007}
\end{figure}

We finally note that the situation is qualitatively similar for 
other values of $\omega$ on the unstable branch. In our last example, we 
choose $\omega=0.007<\omega_{\rm crit}$ within the initial $u_0$ given by \eqref{ini2}. 
We again find that perturbations with an initial mass smaller than $M(Q_{\omega=0.007})$ 
are purely dispersive. Perturbations with mass larger than $M(Q_{\omega=0.007})$ lead to 
damped oscillations around 
some final state, see Fig.~\ref{NLS35_d3solom001m007}. 
The asymptotic final state appears 
to be close to $Q_{\omega=0.044}$ on the stable branch. 
The mass of the unperturbed initial data 
$M(Q_{\omega=0.007})\approx 90.57$ is seen to be bigger than 
$M(Q_{\omega=0.044})\approx 75.37$, showing again 
that a non-negligible part of the initial mass has been radiated away. 

\medskip

The numerical results within this subsection can then be summarized as follows: 
\begin{conjecture} For $\omega<\omega_{c}$, consider initial data of the form
 \[
 u_{0}(x)=Q_{\omega}(x)+\epsilon(|x|), \ \text{ 
  with $\|\epsilon\|_{H^{1}}\ll 1$.}
  \]
 \begin{itemize}
\item[(i)] If $M(u_0)<M(Q_\omega)$, then the solution $u$ to \eqref{eq:nls} is purely dispersive;
\item[(ii)] If $M(u_0)>M(Q_\omega)$, then the solution to
  \eqref{eq:nls} converges, as $t\to +\infty$, to a solitary wave
  $\phi_{\underline \omega} (t,x)= e^{i\underline \omega t}Q_{\underline
    \omega}(x)$ plus radiation, where $Q_{\underline \omega}$ 
is a stable ground state with mass smaller that the unstable
one,
$M(Q_{\underline \omega})<M(Q_\omega)$. 
\end{itemize}
\end{conjecture}

\begin{remark}
The same instability scenario was found (numerically) for perturbed solitary wave solutions to 
the generalized BBM equation in \cite{BMR}, and for a version of NLS with derivative nonlinearity in \cite{AKS}. In particular, analogously to our situation, 
a perturbation which lowered the 
mass of the initial data below the one of the (unstable) solitary wave always resulted in purely dispersive solutions. 
In all of these cases, it 
remains an interesting open question to find a possible selection criterion for the specific value $\underline \omega$ which describes the (stable) asymptotic state $\phi_{\underline \omega}$.
\end{remark}


%

\end{document}